\documentclass[a4paper,reqno]{amsart}
\usepackage{amssymb,amsmath,setspace,centernot}
\usepackage{hyperref}
\newcommand{\Z}{\mathbb{Z}}
\newcommand{\Q}{\mathbb{Q}}
\newcommand{\R}{\mathbb{R}}
\newcommand{\C}{\mathbb{C}}
\newcommand{\N}{\mathbb{N}}

\newcommand\sign{\operatorname{sgn}}
\newcommand\im{\operatorname{Im}}
\newcommand\re{\operatorname{Re}}

\newcommand\smod[1]{\operatorname{mod} #1}
\newcommand\trans{\mathrm{t}}

\let\oldtheta=\theta
\renewcommand\theta{\vartheta}
\newcommand{\modtheta}{\skew{0}\widehat{\vphantom{\rule{1pt}{8.4pt}}\smash{\widehat{\Theta}}}} 
\newcommand\ko{\hspace{.1em}}
\newcommand\steilt{\hspace{1pt}|\hspace{1pt}}

\newtheorem{theorem}{Theorem}[section]
\newtheorem{lemma}[theorem]{Lemma}

\theoremstyle{definition}
\newtheorem{definition}[theorem]{Definition}
\newtheorem{remark}[theorem]{Remark} 

\newenvironment{proofof}[1]{\begin{proof}[Proof of {#1}]}{\end{proof}}

\makeatletter
\@namedef{subjclassname@2020}{\textup{2020} Mathematics Subject Classification}
\makeatother

\numberwithin{equation}{section}

\title[Indefinite Theta Series with (Spherical) Polynomials II]{Theta Series for Quadratic Forms of Signature $(n-1,1)$\\ with (Spherical) Polynomials II}
\date{\today}
\author[C. Roehrig, S. Zwegers]{Christina Roehrig and Sander Zwegers}
\address{Department of Mathematics and Computer Science, University of Cologne, Weyertal 86--90, 50931 Cologne, Germany}
\email{croehrig@math.uni-koeln.de}
\email{szwegers@uni-koeln.de}

\begin{document}

\begin{abstract}
We generalize the construction from \cite{RZ} of theta series for quadratic forms of signature $(n-1,1)$ with homogeneous and spherical polynomials.
Namely, we allow that the parameters $c_1,c_2$, which define the theta series and ensure the convergence of the defining series, are located on the boundary of the cone $C_Q$. 
This enables us to study several interesting examples such as Eisenstein series, modular forms on $\Gamma_0(4)$ which appear during the investigation of quadratic polynomials of a fixed discriminant, and a mock theta function of order 2 that is connected to the generating function of the Hurwitz class numbers $H(8n+7)$.
\end{abstract}

\maketitle

\section{Introduction}
Theta series for positive definite quadratic forms $Q:\R^n\longrightarrow\R$ associated to a function $f:\R^n\longrightarrow\C$, i.\,e. series of the form
\[\Theta_{Q,f}(\tau)=\sum_{\ell\in\Z^n}f(\ell)\,q^{Q(\ell)}\quad(q=e^{2\pi i \tau},\,\im (\tau)>0)\]
play an important role in the construction of modular forms. In particular, if $f$ a spherical polynomial of degree $d$, we obtain a holomorphic modular form of weight $n/2+d$ on some subgroup of $\operatorname{SL}_2(\Z)$ and with some character (see \cite{schoeneberg, ogg, shimura}). For indefinite quadratic forms one has to ensure the convergence of the series $\Theta_{Q,f}(\tau)=\sum_{\ell\in\Z^n}f(\ell)\,q^{Q(\ell)}$. One way to do this is to include majorants as was done by Siegel \cite{siegel}. Further, Vign\'eras \cite{vigneras1,vigneras2} gave a general construction for indefinite theta functions. However, note that both constructions give modular forms that are in general non-holomorphic.

Modular forms also arise from different contexts. Zagier \cite{zagier2} studied sums of powers of quadratic polynomials with integer coefficients and discovered for an even integer $k$ and a fixed real number $x$ a modular form of weight $k+1/2$ on $\Gamma_0(4)$ that has the following expansion:
\[T_x(\tau)=\sum_{\substack{(a,b,c)\in\Z^3\\ax^2+bx+c>0>a}} (ax^2+bx+c)^{k-1}\,q^{b^2-4ac} -\frac{1}{2k}\sum_{m=-\infty}^\infty \overline{B}_k(mx)\,q^{m^2}+\delta_{k,2}\frac{ \kappa(x)}{2}\sum_{m=1}^\infty m^2\,q^{m^2}\]
(See Section \ref{section_quadratic} for the precise definitions of $\overline{B}_k$ and $\kappa$.) The first sum has the form of a theta series, where $(ax^2+bx+c)^{k-1}$ is a spherical polynomial with respect to the binary quadratic form $b^2-4ac$ of signature $(2,1)$. In contrast to the definition for positive definite quadratic forms, here one restricts the summation to a cone in $\Z^3$ determined by $ax^2+bx+c>0>a$. Zagier remarks that this observation ``suggests that there may be an arithmetic theory of theta series attached to indefinite quadratic forms in which the summation runs over the intersection of the lattice with some simplicial cone on which the quadratic form is positive and the result is still a modular form of the expected weight and level''.

Later, this theory was developed by G\"ottsche and Zagier \cite{GZ} and the second author \cite{zwegers} for quadratic forms of signature $(n-1,1)$. In these constructions it is crucial to employ vectors $c_1,c_2\in \R^n$ with $Q(c_i)\leq 0$ and $B(c_1,c_2)<0$, which restrict the summation to a part of the lattice where the quadratic form can be bounded by a positive definite quadratic form, which ensures convergence. 

However, these constructions do not include spherical polynomials of higher degree, so it is not possible to recover, for example, the function $T_x$ in the scope of this theory. In a previous project \cite{RZ}, we extended the definition of the theta functions from \cite{zwegers} to include homogeneous and spherical polynomials.
In \cite{RZ} we focused on the case $Q(c)<0$, constructed a holomorphic theta series which is not modular and a corresponding non-holomorphic modular theta series.
Further, we gave a criterion to determine when these two versions coincide in order to construct holomorphic and almost holomorphic modular forms. 

Since there are many interesting modular forms, which can be realized as theta series with spherical polynomials, but require $c\in \R^n$ to be located on the boundary of the cone, the aim of the present paper is to extend the results from \cite{RZ} to include vectors $c$ with $Q(c)=0$.
Most of the reasoning to establish the absolute convergence of the theta series and its modularity properties are quite analogous to the reasoning in \cite{RZ} and \cite{zwegers}.
More interestingly, we show in Section \ref{section_examples} how we can use the results for $Q(c)=0$ to establish the (mock) modularity in certain special cases: we consider the usual Eisenstein series, embed the aforementioned function $T_x$ in this theory, and show that the generating function of the Hurwitz class numbers $H(8n+7)$ is a mock theta function.

\section{Definitions and statement of the main results}
For the rest of the paper we assume that the quadratic form $Q$ has signature $(n-1,1)$.
We let $A$ denote the corresponding symmetric matrix (so $Q(v)=\frac12 v^\trans Av$), where we assume $A\in \Z^{n\times n}$. Further, let $B$ be the bilinear form associated to $Q$: $B(u,v)=u^\trans Av=Q(u+v)-Q(u)-Q(v)$.
Since $Q$ has signature $(n-1,1)$, the set of vectors $c\in\R^n$ with $Q(c)<0$ has two components.
If $B(c_1,c_2)<0$, then $c_1$ and $c_2$ belong to the same component, while if $B(c_1,c_2)>0$ then $c_1$ and $c_2$ belong to opposite components.
Let $C_Q$ be one of those components.
If $c_0$ is in that component, then $C_Q$ is given by:
\[C_Q :=\{ c\in \R^n \mid Q(c)<0,\ B(c,c_0)<0\}\]
Here we also consider the corresponding set of cusps
\[S_Q:= \{ c\in \Q^n \mid Q(c)=0,\ B(c,c_0)<0\}\]
and let $\overline{C}_Q:=C_Q\cup S_Q$.
For $c\in \overline{C}_Q$ we set
\[R(c):=\begin{cases}
\R^n&\text{if }c\in C_Q,\\
\lbrace a\in \R^n\mid B(c,a)\notin \Z\rbrace &\text{if }c\in S_Q.\\
\end{cases}\]
In \cite{zwegers}, the second author used the error function
\begin{align*}
E(z) := 2\int_0^z e^{-\pi u^2}du=\operatorname{sgn}(z) - \operatorname{sgn}(z) \int_{z^2}^\infty u^{-1/2} e^{-\pi u} du
\end{align*}
to define a non-holomorphic modular theta series and determine its holomorphic part. We generalize this construction as follows (note that this is a slightly more general definition as in \cite{RZ} since we do not necessarily normalize $c\in C_Q$ and include $c\in S_Q$):
\begin{definition}
Let $\Delta=\Delta_Q:=\bigl(\frac{\partial}{\partial v}\bigr)^\trans A^{-1}\frac{\partial}{\partial v}$ denote the Laplacian associated to $Q$ (we often omit $Q$ in the notation, as we take it to be fixed).
We set
\[ e^{-\Delta/8\pi}:= \sum_{k=0}^\infty \frac{(-1)^k}{(8\pi)^k k!}\ko \Delta^k,\qquad \partial_c:= \frac{1}{\sqrt{-Q(c)}}\,c^\trans \frac{\partial}{\partial v}=\frac{1}{\sqrt{-Q(c)}}\,\sum_{i=1}^n c_i \frac{\partial}{\partial v_i}\]
and for a homogeneous polynomial $f:\R^n\longrightarrow\C$ of degree $d$ we define $\widehat f:=e^{-\Delta /8\pi} f$. Further, we set
\[p^{c}[f](v):=\begin{cases}\sum_{k=0}^d \frac{(-1)^k}{(4\pi)^kk!}\ko E^{(k)}\Bigl(\frac{B(c,v)}{\sqrt{-Q(c)}}\Bigr)\cdot \partial_c^k \widehat f(v)&\text{if }c\in C_Q,\\
\sign(B(c,v))\cdot \widehat{f}(v)&\text{if }c\in S_Q.
\end{cases}\]
\end{definition}

\begin{definition}
Let $f:\R^n \longrightarrow \C$ be a homogeneous polynomial of degree $d$ and let $c_1,c_2\in\overline{C}_Q$.
We define the holomorphic theta series associated to $Q$ and $f$ with characteristics $a\in R(c_1)\cap R(c_2)$ and $b\in \R^n$ by
\[\Theta_{a,b}^{c_1,c_2}[f] (\tau):= \sum_{\ell\in a+\Z^n}\bigl\{\sign (B(c_1,\ell)) - \sign (B(c_2,\ell))\bigr\}\ko f(\ell)\ko q^{Q(\ell)}\ko e^{2\pi i B(\ell,b)},\]
the almost holomorphic theta series by
\[\widehat\Theta_{a,b}^{c_1,c_2}[f] (\tau):= y^{-d/2}\sum_{\ell\in a+\Z^n}\bigl\{\sign (B(c_1,\ell)) - \sign (B(c_2,\ell))\bigr\}\ko \widehat{f}(\ell y^{1/2})\ko q^{Q(\ell)}\ko e^{2\pi i B(\ell,b)},\]
and the non-holomorphic theta series by
\[ \modtheta_{a,b}^{c_1,c_2}[f](\tau) := y^{-d/2} \sum_{\ell\in a+\Z^n} \bigl\{p^{c_1}[f](\ell y^{1/2}) - p^{c_2}[f](\ell y^{1/2})\bigr\}\ko q^{Q(\ell)}\ko e^{2\pi iB(\ell,b)}.\]
\end{definition}

\begin{remark}
(a) We show in Lemma \ref{growth} that all three theta series $\Theta_{a,b}^{c_1,c_2}[f],\widehat\Theta_{a,b}^{c_1,c_2}[f],\modtheta_{a,b}^{c_1,c_2}[f]$ are absolutely convergent.
In \cite{RZ} this was already shown for the case $c_1,c_2\in C_Q$.\\
(b) In \cite{RZ} we derived certain conditions under which the non-holomorphic theta series agrees with the almost holomorphic theta series.
For $c_1,c_2\in S_Q$ it immediately follows from the definition that $\modtheta_{a,b}^{c_1,c_2}$ agrees with $\widehat\Theta_{a,b}^{c_1,c_2}$, which is an almost holomorphic theta series of depth $\leq d/2$ since $f$ is a homogeneous polynomial of degree $d$.\\
(c) As usual we call a polynomial $f:\R^n\longrightarrow\C$ spherical (of degree $d$) if it is homogeneous (of degree $d$) and vanishes under the Laplacian, i.\,e. $\Delta f=0$.
If $f$ is spherical of degree $d$, we have $\widehat f=e^{-\Delta /8\pi} f=f$ and $y^{-d/2}\widehat f(\ell y^{1/2})= y^{-d/2} f(\ell y^{1/2})=f(\ell)$.
Hence in this case $\widehat\Theta_{a,b}^{c_1,c_2}[f]$ and the holomorphic theta series $\Theta_{a,b}^{c_1,c_2}[f]$ agree.
\end{remark}

In the next section we'll show that the following two theorems hold for the modular theta series $\modtheta_{a,b}^{c_1,c_2}$.

\begin{theorem}\label{limit}
Let $c_1,c_3\in C_Q$, $c_2\in S_Q$ and $a\in R(c_2)$, $b\in \R^n$. For $c(t)=c_2+tc_3$, we have $c(t)\in C_Q$ for all $t\in (0,\infty)$ and $\lim_{t\downarrow 0} \modtheta_{a,b}^{c_1,c(t)} =\modtheta_{a,b}^{c_1,c_2}$.
\end{theorem}

In the next theorem, we collect the (modular) 
properties of $\modtheta_{a,b}^{c_1,c_2}$:

\begin{theorem}\label{modularity}
Let $c_1,c_2\in \overline{C}_Q$ and $a\in R(c_1),b\in \R^n$. Further, let $A\in \Z^{n\times n}$ be a symmetric matrix.
The theta function with respect to $(c_1,c_2)$ satisfies
\begin{enumerate}
\item $\modtheta_{a+\lambda,b}^{c_1,c_2}[f](\tau) =\modtheta_{a,b}^{c_1,c_2}[f](\tau)\quad\text{for }\lambda\in \Z^n,$
\item $\modtheta_{a,b+\mu}^{c_1,c_2}[f](\tau) =e^{2\pi i B(a,\mu)}\ko\modtheta_{a,b}^{c_1,c_2}[f](\tau)\quad\text{for }\mu\in A^{-1}\Z^n,$
\item $\modtheta_{-a,-b}^{c_1,c_2}[f](\tau) =(-1)^{d+1}\modtheta_{a,b}^{c_1,c_2}[f](\tau),$
\item $\modtheta_{a,b}^{c_1,c_2}[f](\tau +1) = e^{-2\pi iQ(a)-\pi i B(A^{-1}A^{\ast},a)}\ko \modtheta_{a,a+b+\frac12 A^{-1}A^{\ast}}^{c_1,c_2}[f](\tau)$ with $A^{\ast}$ the vector of diagonal elements of $A$, and $\modtheta_{a,b}^{c_1,c_2}[f](\tau+1) = e^{-2\pi iQ(a)}\ko \modtheta_{a,a+b}^{c_1,c_2}[f](\tau)\quad\text{if }A\text{ is even},$
\item For $a,b\in R(c_1)\cap R(c_2)$ we have
\[\modtheta_{a,b}^{c_1,c_2}[f]\bigl(-\frac{1}{\tau}\bigr) = (-i\tau)^{n/2+d}\ko \frac{i^{d+1}}{\sqrt{|\det A|}}\, e^{2\pi i B(a,b)} \sum_{p\in A^{-1}\Z^n\smod \Z^n} \modtheta_{b+p,-a}^{c_1,c_2}[f](\tau).\]
\end{enumerate}
\end{theorem}

\section{Convergence of the theta series and proof of the main results}

We first show the convergence of the theta series.
\begin{lemma}\label{growth}
For $c_1,c_2\in \overline{C}_Q$ the series defining the theta functions $\Theta_{a,b}^{c_1,c_2}[f]$, $\widehat\Theta_{a,b}^{c_1,c_2}[f]$ and $\modtheta_{a,b}^{c_1,c_2}[f]$ are absolutely convergent.
\end{lemma}

\begin{proof}
In Remark 2.4(a) and Lemma 3.1 in \cite{RZ} this is shown for the case $c_1,c_2\in C_Q$.
Note that we slightly changed the definition of the theta series, as we introduced the characteristics $a,b$. However, this does not change the convergence properties.
Using the same argumentation we find that it suffices to show that the series
\begin{align}\label{align_conv}
    \sum_{\ell\in a+\Z^n} \lbrace \sign \bigl(B(c_1,\ell)\bigr)-\sign \bigl(B(c_2,\ell)\bigr)\rbrace\ko P(\ell)\ko q^{Q(\ell)}\ko e^{2\pi i B(\ell,b)}
\end{align}
is absolutely convergent for any polynomial $P$ and $c_1,c_2\in \overline{C}_Q$.
In \cite{RZ} we have already treated the case $c_1,c_2\in C_Q$.
To obtain the other cases, we note that by the cocycle condition $\modtheta_{a,b}^{c_1,c_2}+\modtheta_{a,b}^{c_2,c_3}+\modtheta_{a,b}^{c_3,c_1}=0$ for $c_1,c_2,c_3\in \overline{C}_Q$ it suffices to consider the case $c_1\in C_Q$ and $c_2\in S_Q$ (the claim then also follows for $c_1\in S_Q,\,c_2\in C_Q$ and $c_1,c_2\in S_Q$).
For this, we use more or less the same arguments as in Case 2 of the proof of Proposition 2.4 in \cite{zwegers}: first of all, we can assume that $c_1\in \Z^n\cap C_Q$ and $c_2\in \Z^n\cap S_Q$. We also choose the same decomposition $\ell=\mu +mc_2$ with $\mu \in a+\Z^n$ and $m\in \Z$, such that $\frac{B(c_1,\mu)}{B(c_1,c_2)}\in [0,1)$. Then we can write the series in \eqref{align_conv} as 
\[ -\sum_{\substack{\mu\in a+\Z^n\\ \frac{B(c_1,\mu)}{B(c_1,c_2)}\in [0,1)}}\sum_{m\in \Z} \Big\lbrace \sign \bigl(B(c_2,\mu)\bigr)+\sign \Bigl(m+\frac{B(c_1,\mu)}{B(c_1,c_2)}\Bigr)\Big\rbrace\ko P(\mu +mc_2)\ko q^{Q(\mu)+B(c_2,\mu)m}\ko e^{2\pi i B(\mu,b)+2\pi i B(c_2,b)m}.\]
We have $P(\mu+mc_2)=\sum_{k=0}^d P_k(\mu)\ko m^k$ for some polynomials $P_k:\R^n\longrightarrow \C$.
It suffices to show that the series above converges absolutely for one of these summands with a fixed $k\in\lbrace 0,\ldots, d\rbrace$.
Thus we consider the inner series
\begin{align}\label{align_innersum}
\sum_{m\in \Z} \Big\lbrace \sign \bigl(B(c_2,\mu)\bigr)+\sign \Bigl(m+\frac{B(c_1,\mu)}{B(c_1,c_2)}\Bigr)\Big\rbrace\ko m^k\ko q^{B(c_2,\mu)m}\ko e^{2\pi i B(c_2,b)m}.
\end{align}
For $B(c_2,\mu)>0$, this equation has the form
\[2\sum_{m\geq 0} m^k \ko x^m - \delta_{k,0}\delta_{B(c_1,\mu),0}\quad\text{with }|x|=|q^{B(c_2,\mu)}\ko e^{2\pi i B(c_2,b)}|=e^{-2\pi B(c_2,\mu)y}<1,\]
and for $B(c_2,\mu)<0$, equation \eqref{align_innersum} has the form
\[-2\sum_{m\leq -1} m^k\ko x^m\quad\text{with }|x|=|q^{B(c_2,\mu)}\ko e^{2\pi i B(c_2,b)}|=e^{-2\pi B(c_2,\mu)y}>1.\]
Let $E_{k,j}$ be the Eulerian number, then we have the identities
\begin{align*}
    \frac{x}{(1-x)^{k+1}}\ko \sum_{j=0}^{k-1}E_{k,j}\ko x^j=\begin{cases}
    \sum_{m\geq 0} m^k\ko x^m&\text{if }|x|<1,\\
    -\sum_{m\leq -1} m^k\ko x^m&\text{if }|x|>1.
    \end{cases}
\end{align*}
For $k=0$ this is just the usual geometric series. The second identity follows from the first one for $x^{-1}$ using the fact that $E_{k,j}=E_{k,k-1-j}$ holds.
Thus, equation \eqref{align_innersum} has the form
\[\frac{x}{(1-x)^{k+1}}\ko \sum_{j=0}^{k-1}E_{k,j}\ko x^j\quad\text{with }x=q^{B(c_2,\mu)}\ko e^{2\pi i B(c_2,b)}.\]
The expression $B(c_2,\mu)$ does not become arbitrarily small since $a\in R(c_2)$, so the term above is bounded.
Then one can proceed exactly as in \cite{zwegers} to conclude the proof.
\end{proof}

Now we can prove the two main theorems, which generalize the results in \cite{RZ} to the case $c\in \overline{C}_Q$.

\begin{proofof}{Theorem \ref{limit}}
Again, we can make use of the results on theta series, for which we do not include polynomials in the definition.
We proceed as in \cite{zwegers}:
from Proposition 2.7(5) we know that $c(t)\in C_Q$ for all $t\in (0,\infty)$ and that for the second part of the theorem it is sufficient to show $\lim_{t\downarrow 0}\modtheta_{a,b}^{c_2,c(t)}= 0$ for $c(t)=c_2+tc_1$. For $\bigl(p^{c_2}[f]-p^{c(t)}[f]\bigr)(v)$ we have the decomposition
\begin{align}\label{decomp}
    \bigl(\sign(B(c_2,v))-\sign(B(c(t),v))\bigr)\ko \widehat{f}(v)+\sign(B(c(t),v))\ko \beta\Bigl(\frac{B(c(t),v)^2}{-Q(c(t))}\Bigr)\ko  \widehat{f}(v)-\widetilde{p}^{c(t)}[f](v),
\end{align}
where $\widetilde{p}^{c}[f](v):=\sum_{k=1}^d \frac{(-1)^k}{(4\pi)^kk!}\ko E^{(k)}\Bigl(\frac{B(c,v)}{\sqrt{-Q(c)}}\Bigr)\cdot \partial_c^k \widehat f(v)$ for $c\in C_Q$.

The first two summands are the terms that were also considered in \cite{zwegers} except for the polynomial factor $\widehat{f}$.
However, since $\widehat{f}$ is independent of $t$, we can use the same argumentation as in \cite{zwegers} and use Lemma \ref{growth} to get
\[\lim_{t\downarrow 0}\sum_{\ell\in a+\Z^n}\big\lbrace\sign(B(c_2,\ell))-\sign(B(c(t),\ell))\big\rbrace\ko \widehat{f}(\ell)\ko q^{Q(\ell)}e^{2\pi i B(\ell,b)}=0,\]
and
\[\lim_{t\downarrow 0}\sum_{\ell\in a+\Z^n}\sign(B(c(t),\ell))\ko \beta\Bigl(\frac{B(c(t),\ell)^2}{-Q(c(t))}\Bigr)\ko  \widehat{f}(\ell)\ko q^{Q(\ell)}e^{2\pi i B(\ell,b)}=0,\]
since both series are uniformly convergent for $t\in (0,\infty)$.

As $c(t)\in C_Q$ for $t\in (0,\infty)$, we observe (exactly as in the proof of Lemma 3.1 in \cite{RZ}) that $\widetilde{p}^{c(t)}[f](v)$ can be written as a polynomial in $c(t)/\sqrt{-Q(c(t))}$ and $v$ times the non-polynomial factor $e^{\pi \frac{B(c(t),v)^2}{Q(c(t))}}$. Hence, we have to show that for any polynomial $P$ there exists a majorant for
\[P\Bigl(\frac{c(t)}{\sqrt{-Q(c(t))}},v\Bigr)\ko e^{-2\pi Q(v)+\pi \frac{B(c(t),v)^2}{Q(c(t))}}\]
that is independent of $t$ and for which the sum is absolutely convergent. We have to consider what happens if $t\downarrow 0$, so we now assume $t\in (0,t_0)$ for some $t_0>0$. We follow the proof of Proposition 2.7(5) in \cite{zwegers} and use the same decomposition of $a+\Z^n$ into the three subsets
\begin{align*}
\begin{split}
    P_1&:=\{v\in a+\Z^n\mid \sign(B(c_2,v))=-\sign(B(c_1,v))\},\\
    P_2&:=\{v\in a+\Z^n\mid B(c_1,v)\ko \bigl(B(c_1,c_2) B(c_1,v)-2Q(c_1)B(c_2,v)\bigr)\geq 0\},\\
    P_3&:=\{v\in a+\Z^n\mid \sign(B(c_2,v))=-\sign\bigl(B(c_1,c_2)B(c_1,v)-2Q(c_1)B(c_2,v)\bigr)\}.
\end{split}
\end{align*}
We then determine a majorant on each subset $P_i$ separately.
We show that there exists a polynomial $\widetilde{P}$ such that
\begin{align}\label{eq:Ptilde}
    \bigg|P\Bigl(\frac{c(t)}{\sqrt{-Q(c(t))}},v\Bigr)\bigg|\leq \widetilde{P}\Bigl(\Big|\frac{B(c(t),v)}{\sqrt{-Q(c(t))}}\Big|,|v|\Bigr)
\end{align}
holds for all $v\in \R^n$ ($|v|$ stands for $|v_1|,\ldots,|v_n|$ here) and $t\in (0,t_0)$: We use induction on the degree of $P$ as a polynomial in $c(t)$ and first note that we can write $P$ as a finite linear combination of terms of the form $B\big(\xi,\frac{c(t)}{\sqrt{-Q(c(t))}}\big)\ko P_{\xi}\Bigl(\frac{B(c(t),v)}{\sqrt{-Q(c(t))}},v\Bigr)$ for suitable $\xi\in\C$, where the degree of $P_\xi$ is strictly lower than the degree of $P$ in $c(t)$. So we can use the induction hypothesis on $P_\xi$, and we also observe
\[\begin{split}\bigg|B\Big(\xi,\frac{c(t)}{\sqrt{-Q(c(t))}}\Big)\bigg|&\leq\Big| \frac{B(\xi,c_2)}{B(c_2,v)}\ko\frac{B(c(t),v)}{\sqrt{-Q(c(t))}}\Big| +\frac{t}{\sqrt{-Q(c(t))}}\ko\bigg| \frac{B(\xi,c_1)B(c_2,v)-B(\xi,c_2)B(c_1,v)}{B(c_2,v)} \bigg|\\
&\leq\alpha\Big|\frac{B(c(t),v)}{\sqrt{-Q(c(t))}}\Big| + R(|v|)
\end{split}\]
for a constant $\alpha>0$ and some polynomial $R$. In the last step, we use that $|B(c_2,v)|$ does not become arbitrarily small for $v\in a+\Z^n$ since we assume $a\in R(c_2)$, and that $t/\sqrt{-Q(c(t))}\rightarrow 0$ for $t\downarrow 0$, so that this estimate holds for $t\in (0,t_0)$.

Using \eqref{eq:Ptilde}, we can now consider
\[\widetilde{P}\Bigl(\Big|\frac{B(c(t),v)}{\sqrt{-Q(c(t))}}\Big|,|v|\Bigr)\ko e^{\pi \frac{B(c(t),v)^2}{Q(c(t))}}.\]
This expression has polynomial growth in $v$ for $t\downarrow 0$ since it is clear that $\Big|\frac{B(c(t),v)}{\sqrt{-Q(c(t))}}\Big|^{\alpha}\ko e^{\pi \frac{B(c(t),v)^2}{Q(c(t))}}$ is bounded for any $\alpha\in \N_0$.
So we find a polynomial $S$ as an upper bound that is independent of $t$. In Lemma \ref{growth} we have seen that the sum
\[\sum_{v\in P_1} S(v)\ko e^{-2\pi Q(v)}\]
is absolutely convergent.

On $P_2$ we observe that
\[\frac{B(c(t),v)^2}{Q(c(t))}\leq \frac{B(c_1,v)^2}{Q(c_1)}+\frac{B(c_2,v)^2}{Q(c(t))}\]  holds.
Thus we have
\[\bigg|P\Bigl(\frac{c(t)}{\sqrt{-Q(c(t))}},v\Bigr)\bigg|\ko e^{-2\pi Q(v)+\pi \frac{B(c(t),v)^2}{Q(c(t))}}\leq \bigg|P\Bigl(\frac{c(t)}{\sqrt{-Q(c(t))}},v\Bigr)\bigg|\ko e^{\pi \frac{B(c_2,v)^2}{2Q(c(t))}}\ko e^{-2\pi \bigl(Q(v)- \frac{B(c_1,v)^2}{2Q(c_1)}\bigr)}.\]
Similarly as above, we have
\[\bigg|P\Bigl(\frac{c(t)}{\sqrt{-Q(c(t))}},v\Bigr)\bigg|\leq \widetilde{P}\Bigl(\Big|\frac{B(c_2,v)}{\sqrt{-Q(c(t))}}\Big|,|v|\Bigr)\]
for some polynomial $\widetilde{P}$ and conclude that
\[ \widetilde{P}\Bigl(\Big|\frac{B(c_2,v)}{\sqrt{-Q(c(t))}}\Big|,|v|\Bigr)\ko e^{\pi \frac{B(c_2,v)^2}{2Q(c(t))}}\] has polynomial growth in $v$ for $t\downarrow 0$. Since the quadratic form $v\mapsto Q(v)- \frac{B(c_1,v)^2}{2Q(c_1)}$ is positive definite (see Lemma 2.5 in \cite{zwegers}), we have constructed a suitable majorant on $P_2$.

On $P_3$ we consider the quadratic form $\widetilde{Q}$ of signature $(n-1,1)$ that is defined in \cite{zwegers} as follows:
\[\widetilde{Q}(v):=Q(v)-\frac{2B(c_2,v)}{B(c_1,c_2)^2}\ko \bigl(B(c_1,c_2)B(c_1,v)-Q(c_1)B(c_2,v)\bigr)\]
We denote by $\widetilde{B}$ the associated bilinear form.
Setting \[\widetilde{c}(\ko\widetilde{t}\ko)=\widetilde{c}_2+\widetilde{t}\ko\widetilde{c}_1\qquad\text{with}\quad\widetilde{c}_1=\frac{B(c_1,c_2)}{2Q(c_1)}c_1-c_2,\quad\widetilde{c}_2=-c_2\quad\text{and}\quad \widetilde{t}=\frac{2Q(c_1)}{B(c_1,c_2)}t,\]
we have
\[Q(v)-\frac{B(c(t),v)^2}{2Q(c(t))}=\widetilde{Q}(v)-\frac{\widetilde{B}(\widetilde{c}(\ko\widetilde{t}\ko ),v)^2}{2\widetilde{Q}(\widetilde{c}(\widetilde{t}))}.\]
(This identity follows quite easily using the relations shown in \cite{zwegers}.) Again, we can find a polynomial $\widetilde{P}$ with
\[\bigg|P\Bigl(\frac{c(t)}{\sqrt{-Q(c(t))}},v\Bigr)\bigg|\leq \widetilde{P}\Biggl(\Bigg|\frac{\widetilde{B}(\widetilde{c}(\ko\widetilde{t}\ko),v)}{\sqrt{-\widetilde{Q}(\widetilde{c}(\widetilde{t}))}}\Bigg|,|v|\Biggr)\]
and then conclude analogously as for $P_1$ that it is sufficient to show the absolute convergence of
\[\sum_{v\in P_3} S(v)\ko e^{-2\pi \widetilde{Q}(v)}\] for some polynomial $S$. This follows for $S\equiv 1$ as in \cite{zwegers} and then also for non-constant polynomials by Lemma \ref{growth}.

Due to the uniform convergence in $t\in(0,t_0)$, we can then consider $\lim_{t\downarrow 0}\widetilde{p}^{c(t)}[f](v)$ and see that the convergence is dominated by the part $e^{\pi \frac{B(c(t),v)^2}{Q(c(t))}}$, which goes to zero for $t\downarrow 0$.
We have thus shown 
\[\lim_{t\downarrow 0}\sum_{\ell\in a+\Z^n}\widetilde{p}^{c(t)}[f](\ell)\ko q^{Q(\ell)}e^{2\pi i B(\ell,b)}=0.\]
Combining the results for the three separate series that we have obtained by the decomposition \eqref{decomp}, we conclude that $\lim_{t\downarrow 0}\modtheta_{a,b}^{c_2,c(t)}=0$ holds.
\end{proofof}

To describe the modular transformation behavior of $\modtheta_{a,b}^{c_1,c_2}$ we can employ Vignéras' result \cite{vigneras2}.
For the proof, we consider the Fourier transform $\bigl(\mathcal{F}f\bigr)(v):=\int_{\R^n} f(u)\ko e^{-2\pi i u\cdot v}\ko du$.
\begin{proofof}{Theorem \ref{modularity}}
(1), (2) and (4) immediately follow by the same calculations as in \cite{zwegers}.
The third identity holds since $p^c[f](-v)=(-1)^{d+1} p^c[f](v)$ for $c\in \overline{C}_Q$.
To show part (5) we first consider $c_1,c_2\in C_Q$ and make use of Vignéras' result \cite{vigneras2}:
let $g:=p^{c_1}[f]-p^{c_2}[f]$. From Lemma 3.2 of \cite{RZ} we know that $Dg=d\ko g$, so from \cite{vigneras2} we get that the following identity holds for $g_{\tau}(u):=y^{-d/2}\ko g(uy^{1/2})\ko e^{2\pi iQ(u)\tau}$:
\[\bigl(\mathcal{F}g_{-1/\tau}\bigr)(v)=(-i\tau)^{n/2+d}\ko \frac{(-i)^{d+1}}{\sqrt{|\det A|}}\ko g_{\tau}(A^{-1}v)\]
Hence the Fourier transform of $v\mapsto g_{-1/\tau}(v+a)\ko e^{2\pi i B(v+a,b)}$ is
\[ \int_{\R^n} g_{-1/\tau}(u+a)\ko e^{2\pi i B(u+a,b)-2\pi iu\cdot v}\ko du=(-i\tau)^{n/2+d}\ko \frac{(-i)^{d+1}}{\sqrt{|\det A|}}\ko e^{2\pi ia\cdot v} g_{\tau}(A^{-1}v-b)\]
By applying the Poisson summation formula, which states 
\[\sum_{v\in \Z^n} f(v)=\sum_{v\in \Z^n} \bigl(\mathcal{F}f\bigr)(v),\]
we thus obtain:
\[\begin{split}
\modtheta_{a,b}^{c_1,c_2}[f]\bigl(-\frac{1}{\tau}\bigr) &=
\sum_{v\in\Z^n} g_{-1/\tau}(v+a)\ko e^{2\pi iB(v+a,b)}\\
&= (-i\tau)^{n/2+d}\ko \frac{(-i)^{d+1}}{\sqrt{|\det A|}} \sum_{v\in\Z^n}e^{2\pi ia\cdot v}  g_{\tau}(A^{-1}v-b)\\
&= (-i\tau)^{n/2+d}\ko \frac{(-i)^{d+1}}{\sqrt{|\det A|}}\ e^{2\pi iB(a,b)} \sum_{u\in -b+A^{-1}\Z^n}e^{2\pi iB(a,u)}  g_{\tau}(u)\\
&= (-i\tau)^{n/2+d}\ko \frac{(-i)^{d+1}}{\sqrt{|\det A|}}\ e^{2\pi iB(a,b)} \sum_{p\in A^{-1}\Z^n\smod \Z^n}\modtheta_{-b-p,a}^{c_1,c_2}[f](\tau)
\end{split}\]
With part (3) we then get the desired result for the case $c_1,c_2\in C_Q$.
Using Lemma \ref{limit} we get that the formula also holds for the case $c_1\in C_Q$ and $c_2\in S_Q$.
By the cocycle condition the remaining two cases then follow immediately.
\end{proofof}

\section{Examples}\label{section_examples}
The next lemma proves to be helpful in the construction of the following examples.
\begin{lemma}\label{lemma_bernoulli}
Let $k\in \N$ and let $B_k$ denote the usual $k$-th Bernoulli polynomial. For all $\alpha,\beta\in \R$, and $z\in \C$ with $|z|<2\pi$ we have
\[\frac{(-1)^{k}\ko (k-1)!}{z^k}-\sum_{m=0}^\infty \frac{B_{m+k}(\alpha+\beta-\lfloor\beta\rfloor)}{m+k}\frac{z^m}{m!}=\begin{cases}
\sum_{n+\beta\geq 0} (n+\alpha+\beta)^{k-1}\ko e^{(n+\alpha+\beta)z}&\text{if }\re (z)<0,\\
-\sum_{n+\beta\leq -1} (n+\alpha+\beta)^{k-1}\ko e^{(n+\alpha+\beta)z}&\text{if }\re (z)>0.
\end{cases}\]
\end{lemma}
\begin{proof}
We note that
\[(n+\alpha+\beta)^{k-1}\ko e^{(n+\alpha+\beta)z}=\Bigl(\frac{\partial}{\partial z}\Bigr)^{k-1}e^{(n+\alpha+\beta)z}\]
holds.
Well-known identities for the geometric series yield for any $x\in \C$ and $\beta\in \R$
\[\frac{x^{-\lfloor \beta \rfloor}}{1-x}=\begin{cases}
\sum_{n+\beta\geq 0} x^n &\text{if }|x|<1,\\
-\sum_{n+\beta\leq -1} x^n&\text{if }|x|>1.
\end{cases}\]
As these series are uniformly convergent, we can interchange summation and differentiation and obtain for $\re (z)<0$
\begin{align*}
 \sum_{n+\beta\geq 0} (n+\alpha+\beta)^{k-1}\ko e^{(n+\alpha+\beta)z}=\Bigl(\frac{\partial}{\partial z}\Bigr)^{k-1}\Bigl( e^{(\alpha+\beta) z}\sum_{n+\beta\geq 0} e^{nz}\Bigr)=\Bigl(\frac{\partial}{\partial z}\Bigr)^{k-1}\Bigl( \frac{e^{(\alpha+\beta-\lfloor\beta\rfloor) z}}{1-e^z}\Bigr)
\end{align*}
and for $\re (z)>0$
\begin{align*}
  -\sum_{n+\beta\leq -1} (n+\alpha+\beta)^{k-1}\ko e^{(n+\alpha+\beta)z}=\Bigl(\frac{\partial}{\partial z}\Bigr)^{k-1}\Bigl( \frac{e^{(\alpha+\beta-\lfloor\beta\rfloor) z}}{1-e^z}\Bigr).
\end{align*}
For $|z|<2\pi$, we know the generating function for the Bernoulli polynomials $B_m$:
\[\frac{e^{(\alpha+\beta-\lfloor\beta\rfloor) z}}{1-e^z}=-\sum_{m=0}^\infty B_{m}(\alpha+\beta-\lfloor\beta\rfloor)\ko\frac{z^{m-1}}{m!}\]
From
\[\Bigl(\frac{\partial}{\partial z}\Bigr)^{k-1} \Bigl( -\sum_{m=0}^\infty B_{m}(\alpha+\beta-\lfloor\beta\rfloor)\ko\frac{z^{m-1}}{m!}\Bigr)=\frac{(-1)^{k}\ko (k-1)!}{z^k}-\sum_{m=k}^\infty \frac{B_{m}(\alpha+\beta-\lfloor\beta\rfloor)}{m}\ko\frac{z^{m-k}}{(m-k)!}\]
the claim follows by shifting $m$ to $m+k$.
\end{proof}
\subsection{Eisenstein Series}
For the Eisenstein series of positive even weight $k$ we use the normalized version
\[G_k(\tau)= -\frac{B_k}{2k}+\sum_{n\in \N^2} n_1^{k-1}q^{n_1n_2}\]
(with $B_k$ the $k$-th Bernoulli number) as definition.

The matrix $A=\bigl(\begin{smallmatrix}0&1\\1&0\end{smallmatrix}\bigr)$ of signature $(1,1)$ induces the quadratic form $Q(v)=v_1v_2$ and the associated bilinear form $B(u,v)=u_1v_2+u_2v_1$. For $c_1=\bigl(\begin{smallmatrix}0\\1\end{smallmatrix}\bigr)$ and $c_2=\bigl(\begin{smallmatrix}-1\\0\end{smallmatrix}\bigr)$ both in $S_Q$, $a,b\in R(c_1)\cap R(c_2)$ and $f$ the spherical polynomial $f(v)=v_1^{k-1}$ we obtain the holomorphic theta series
\begin{align*}
    \Theta_{a,b}^{c_1,c_2}(\tau)=\modtheta_{a,b}^{c_1,c_2}(\tau)=\sum_{\ell \in a+\Z^2}\lbrace\sign(\ell_1)+\sign(\ell_2)\rbrace\ko \ell_1^{k-1}\ko q^{\ell_1 \ell_2}\ko e^{2\pi i (\ell_1 b_2+\ell_2b_1)}.
\end{align*}
By Theorem \ref{modularity} this theta series transforms as follows ($A$ is an even unimodular matrix):
\begin{align*}
\Theta_{a,b}^{c_1,c_2}(\tau+1)= e^{-2\pi i Q(a)}\ko \Theta_{a,a+b}^{c_1,c_2}(\tau)\quad\text{and}\quad
    \Theta_{a,b}^{c_1,c_2}\Bigl(-\frac{1}{\tau}\Bigr)=\tau^k\ko e^{2\pi i B(a,b)}\ko \Theta_{b,-a}^{c_1,c_2}(\tau)
\end{align*}
Next we consider what happens if we let $a,b\rightarrow 0$ (note that $0\notin R(c_1)\cap R(c_2)$).
\begin{lemma}
The meromorphic function 
\[f_{a,b}^G (\tau):= e^{2\pi i a_2 b_1}\ko \frac{(k-1)!}{(2\pi i (a_2\tau+b_2))^k}\]
has the same modular transformation behavior as $\Theta_{a,b}^{c_1,c_2}$ on $\Gamma_1$.
\end{lemma}
\begin{proof}
We derive the modular transformation behavior by the following two straightforward calculations: we have
\begin{align*}
    f_{a,b}^G (\tau+1)=e^{-2\pi i a_1a_2}\ko e^{2\pi i a_2 (a_1+b_1)}\ko \frac{(k-1)!}{(2\pi i (a_2\tau+a_2+b_2))^k}=e^{-2\pi i Q(a)}\ko f_{a,a+b}^G (\tau)
\end{align*}
and
\[f_{a,b}^G\Bigl(-\frac{1}{\tau}\Bigr)=\tau^k\ko e^{2\pi i (a_1b_2+a_2b_1)}\ko e^{-2\pi i a_1b_2} \ko \frac{(k-1)!}{(2\pi i (b_2\tau-a_2))^k} =\tau^k\ko e^{2\pi i B(a,b)}\ko f_{b,-a}^G(\tau).\qedhere\]
\end{proof}
\begin{lemma}\label{prop_eisenstein}
We have
\begin{align*}
    G_k(\tau)=\lim_{a,b\rightarrow 0} \Bigl( \frac14 \Theta_{a,b}^{c_1,c_2}(\tau)-\frac12 f_{a,b}^G (\tau)\Bigr)\quad(k\geq 4),
\end{align*}
and for $k=2$
\begin{align*}
    G_2(\tau)-\frac{1}{4\pi i \tau}=\lim_{a\rightarrow 0}\Bigl( \lim_{b\rightarrow 0} \Bigl( \frac14 \Theta_{a,b}^{c_1,c_2}(\tau)-\frac12 f_{a,b}^G (\tau)\Bigr)\Bigr),\quad  G_2(\tau)=\lim_{b\rightarrow 0}\Bigl( \lim_{a\rightarrow 0} \Bigl( \frac14 \Theta_{a,b}^{c_1,c_2}(\tau)-\frac12 f_{a,b}^G (\tau)\Bigr)\Bigr).
\end{align*}
\end{lemma}

\begin{proof}
To ensure $a,b\in R(c_1)\cap R(c_2)$, we assume $a_i,b_i\in (-1,0)\cup (0,1)$ for $i=1,2$. Further, we write $a+\Z^2$ as the disjoint union of $L_1,L_2,L_3$ with
\begin{align*}
    L_1:=\big\lbrace a+\bigl(\begin{smallmatrix}\ell_1\\ \ell_2\end{smallmatrix}\bigr)\steilt \ell_1,\ell_2\in \Z\setminus \lbrace 0\rbrace \big\rbrace,\quad
    L_2:=\big\lbrace \bigl(\begin{smallmatrix}a_1+\ell_1\\ a_2\end{smallmatrix}\bigr)\steilt \ell_1\in \Z\setminus\lbrace 0\rbrace\big\rbrace,\quad
    L_3:=\big\lbrace \bigl(\begin{smallmatrix}a_1\\ a_2+\ell_2\end{smallmatrix}\bigr)\steilt \ell_2\in \Z\big\rbrace.
\end{align*}
We consider $\frac14\Theta_{a,b}^{c_1,c_2}$, where we first restrict the summation to $L_1$.
Letting $a,b\rightarrow 0$ we get (for $k\geq 2$ even):
\begin{align*}
   \frac14\sum_{\substack{\ell\in \Z^2\\ \ell_1,\ell_2\neq 0}} \lbrace \sign (\ell_1)+\sign (\ell_2) \rbrace\ko \ell_1^{k-1}\ko q^{\ell_1\ell_2}=\sum_{\ell\in \N^2} \ell_1^{k-1}\ko q^{\ell_1\ell_2}
\end{align*}
Restricting the summation to $L_2$ we obtain the expression 
\begin{align*}
  \frac{e^{2\pi i a_2 b_1}}{4}\sum_{\ell_1\in \Z\setminus\lbrace 0\rbrace} \lbrace \sign (a_1+\ell_1)+\sign (a_2) \rbrace\ko (a_1+\ell_1)^{k-1}\ko q^{(a_1+\ell_1)a_2}\ko e^{2\pi i (a_1+\ell_1)b_2}.
\end{align*}
Using Lemma \ref{lemma_bernoulli} (for even $k$) with $\alpha=a_1$, $\beta=0$ and $z=2\pi i (a_2\tau +b_2)$ this equals
\begin{align*}
   \frac{e^{2\pi i a_2 b_1}}{2}\Bigl( \frac{(k-1)!}{(2\pi i (a_2\tau +b_2))^k}-\sum_{m=0}^\infty \frac{B_{m+k}(a_1)}{m+k}\frac{(2\pi i (a_2\tau +b_2))^m}{m!} \Bigr).
\end{align*}
Note that we have to add the summand for $\ell_1=0$ if $a_2$ is positive in order to apply Lemma \ref{lemma_bernoulli}, but as this extra term goes to zero if $a$ goes to zero, we immediately neglect it here.
We subtract the non-holomorphic part $\frac12 f_{a,b}^G$, and for the remaining part we obtain
\begin{align*}
   \lim_{a,b\rightarrow 0} \Bigl(\frac{-e^{2\pi i a_2 b_1}}{2}\sum_{m=0}^\infty \frac{B_{m+k}(a_1)}{m+k}\frac{(2\pi i (a_2\tau +b_2))^m}{m!} \Bigr)=-\frac{B_k}{2k}.
\end{align*}
On $L_3$ we just use the usual identity for the geometric series and thus have
\begin{multline}\label{align_limit}
  \lim_{a,b\rightarrow 0} \Bigl( \frac{1}{4} q^{a_1a_2}e^{2\pi i (a_1b_2+a_2b_1)} a_1^{k-1}\sum_{\ell_2\in \Z} \lbrace \sign (a_1)+\sign (a_2+\ell_2) \rbrace\ko e^{2\pi i (a_1\tau+b_1)\ell_2}\Bigr)\\
  =\lim_{a,b\rightarrow 0} \Bigl(\frac{1}{2} q^{a_1a_2}e^{2\pi i (a_1b_2+a_2b_1)} a_1^{k-1}\frac{1}{1-q^{a_1} e^{2\pi i b_1}}\Bigr).
\end{multline}
If $a_2$ is negative, we again have to add an extra term for $\ell_2=0$, which also goes to zero if $a$ goes to zero.
If $k\geq 4$, we immediately see that \eqref{align_limit} is zero, as $a_1^{k-1}$ has a zero of order higher or equal to three and we only have a simple pole. If $k=2$, it plays a role in which order we take the limit: we have
\[\lim_{b\rightarrow 0} \Bigl(\lim_{a\rightarrow 0} \Bigl(\frac{1}{2} q^{a_1a_2}e^{2\pi i (a_1b_2+a_2b_1)} a_1\frac{1}{1-q^{a_1} e^{2\pi i b_1}}\Bigr)\Bigr)=0\]
since $1-e^{2\pi i b_1}\neq 0$. On the other hand, 
we have
\[\lim_{a\rightarrow 0} \Bigl(\lim_{b\rightarrow 0} \Bigl(\frac{1}{2} q^{a_1a_2}e^{2\pi i (a_1b_2+a_2b_1)} a_1\frac{1}{1-q^{a_1} e^{2\pi i b_1}}\Bigr)\Bigr)=\frac{1}{2}\,\lim_{a\rightarrow 0} \frac{a_1\ko q^{a_1a_2}}{1-q^{a_1}}=-\frac{1}{4\pi i \tau}\]
by L'Hospital's rule.
Combining the different parts then gives the desired result.
\end{proof}

Using the identities in Lemma \ref{prop_eisenstein} we can now easily recover the well-known modular transformation properties of $G_k$.

\begin{theorem}
For $k\geq 4$ the Eisenstein series $G_k$ is a modular form of weight $k$ on $\Gamma_1$: it satisfies $G_k(\tau+1)=G_k(\tau)$ and $G_k(-1/\tau)=\tau^k G_k(\tau)$. Further, $G_2$ satisfies $G_2(\tau+1)=G_2(\tau)$ and
\[G_2(\tau)-\frac{1}{\tau^2}G_2\Bigl(-\frac{1}{\tau}\Bigr)=\frac{1}{4\pi i \tau}.\]
\end{theorem}

\begin{proof}
 For $k\geq 4$ we have
\begin{align*}
    G_k(\tau +1)=\lim_{a,b\rightarrow 0} \Bigl( \Bigl(\frac14 \Theta_{a,b}^{c_1,c_2}-\frac12 f_{a,b}^G\Bigr)(\tau+1)\Bigr)= \lim_{a,b\rightarrow 0} \Bigl( e^{-2\pi i Q(a)}\ko \Bigl(\frac14 \Theta_{a,a+b}^{c_1,c_2}-\frac12 f_{a,a+b}^G\Bigr) (\tau)\Bigr)= G_k(\tau)
\end{align*}
and
\begin{align*}
    G_k\Bigl(-\frac{1}{\tau}\Bigr)=\lim_{a,b\rightarrow 0} \Bigl( \Bigl(\frac14 \Theta_{a,b}^{c_1,c_2}-\frac12 f_{a,b}^G\Bigr)\Bigl(-\frac{1}{\tau}\Bigr)\Bigr)= \lim_{a,b\rightarrow 0} \Bigl( \tau^k e^{2\pi i B(a,b)}\ko \Bigl(\frac14 \Theta_{b,-a}^{c_1,c_2}-\frac12 f_{b,-a}^G\Bigr) (\tau)\Bigr)=\tau^k\ko G_k(\tau),
\end{align*}
which shows that $G_k$ is a holomorphic modular form of weight $k$ on $\Gamma_1$.

For $k=2$ we get $G_2(\tau+1)=G_2(\tau)$ using the same identity as for $k\geq 4$, when we obey the order of the limits and take $\lim_{b\rightarrow 0} \lim_{a\rightarrow 0}$. Further,
\[\begin{split}
   \frac{1}{\tau^2}G_2\Bigl(-\frac{1}{\tau}\Bigr)&= \frac{1}{\tau^2}\lim_{b\rightarrow 0}\Bigl( \lim_{a\rightarrow 0} \Bigl( \Bigl(\frac14 \Theta_{a,b}^{c_1,c_2}-\frac12 f_{a,b}^G \Bigr)\Bigl(-\frac{1}{\tau}\Bigr)\Bigr)\Bigr)\\
   &= \lim_{b\rightarrow 0}\Bigl(\lim_{a\rightarrow 0} \Bigl(e^{2\pi i B(a,b)}\ko \Bigl(\frac14 \Theta_{b,-a}^{c_1,c_2}-\frac12 f_{b,-a}^G\Bigr) (\tau)\Bigr)\Bigr)= G_2(\tau)-\frac{1}{4\pi i \tau}.\qedhere
\end{split}\]
\end{proof}

\subsection{Quadratic Polynomials}\label{section_quadratic}
In this example, we consider two modular forms of weight $k+1/2$ that were discussed by Zagier in \cite{zagier2}.
One considers the quadratic polynomial $ax^2+bx+c$ with $a,b,c\in \Z$ and discriminant $D:=b^2-4ac$. In the following, $B_k$ denotes the $k$-th Bernoulli polynomial and $\overline{B}_k(x):=B_k(x-\lfloor x \rfloor )$ the periodic version of the Bernoulli polynomial.

For $k\in \N$ even and $D$ not a square, let \[P_{k,D}(x):=\sum_{\substack{b^2-4ac=D\\ a>0>c}} (ax^2+bx+c)^{k-1} \in \Z[x].\]
For $D=m^2>0$, we define the right-hand side as $P^{\ast}_{k,m^2}$ and set
\[P_{k,m^2}(x)=P^{\ast}_{k,m^2}(x) + \frac{1}{k}\Bigl(B_k(mx)-x^{2k-2}B_k\bigl(\frac{m}{x}\bigr)\Bigr).\]
For $D=0$ we simply define
\[ P_{k,0}(x)= (1-x^{2k-2}) \frac{B_k}{2k}.\]
Analogously, we define \[F_{k,D}(x):=\sum_{\substack{b^2-4ac=D\\ ax^2+bx+c>0>a}} (ax^2+bx+c)^{k-1} \in \Z[x]\]
and
\[\begin{split}
F_{k,m^2}(x)&=F^{\ast}_{k,m^2}(x) - \frac{1}{k}\overline{B}_k(mx)+\delta_{k,2}\frac{m^2 \kappa(x)}{2}\quad\text{with }\kappa (x)=\begin{cases}
\frac{1}{s^2}&\text{for }x=\frac{r}{s}\,(\gcd (r,s)=1),\\
0&\text{for }x\in \R\setminus \Q,
\end{cases}\\
F_{k,0}(x)&= -\frac{B_k}{2k}.
\end{split}\]
Summing over all positive discriminants, we consider the generating functions
\begin{align*}
    S_x(\tau):=\sum_{D\geq 0} P_{k,D}(x)\ko q^D \quad\text{and}\quad T_x(\tau):=\sum_{D\geq 0} F_{k,D}(x)\ko q^D\qquad\text(x\in \R)
\end{align*}
and by plugging in the definition of $P_{k,D}$ and $F_{k,D}$ we obtain the expansions
\[\begin{split}
S_x(\tau)&=\sum_{\substack{(a,b,c)\in\Z^3\\a>0>c}} (ax^2+bx+c)^{k-1}\,q^{b^2-4ac} +\frac{1}{2k}\sum_{m=-\infty}^\infty B_k(mx)\,q^{m^2}-x^{2k-2}\frac{1}{2k}\sum_{m=-\infty}^\infty B_k\bigl(\frac mx\bigr)\,q^{m^2},\\
T_x(\tau)&=\sum_{\substack{(a,b,c)\in\Z^3\\ax^2+bx+c>0>a}} (ax^2+bx+c)^{k-1}\,q^{b^2-4ac} -\frac{1}{2k}\sum_{m=-\infty}^\infty \overline{B}_k(mx)\,q^{m^2}+\delta_{k,2}\frac{ \kappa(x)}{2}\sum_{m=1}^\infty m^2\,q^{m^2}.
\end{split}\]
The goal here is to recover the results of Zagier \cite{zagier2} on the modularity of these functions by using these expansions  and the theory of theta series for quadratic forms of signature $(n-1,1)$.
\begin{theorem}\label{mainprop}
We have:
\begin{enumerate}
    \item For $x\in \R$, $S_x$ is a modular form of weight $k+1/2$ on $\Gamma_0(4)$;
    \item For $x\in \Q$, $T_x$ is a modular form of weight $k+1/2$ on $\Gamma_0(4)$.
\end{enumerate}
\end{theorem}
The indefinite theta function that is associated to these two functions is constructed as follows.
The matrix $A=\Bigl(\begin{smallmatrix}0&0&-4\\0&2&0\\-4&0&0\end{smallmatrix}\Bigr)$ defines the quadratic form $Q(n)=n_2^2-4n_1n_3$ of signature $(2,1)$. Further, let $k\in \N$ be even and choose the polynomial $f(n)=(n_1x^2+n_2x+n_3)^{k-1}$, which is spherical of degree $k-1$ with respect to $Q$.
While these parameters stay the same when considering $S_x$ and $T_x$, we choose different elements $c_1,c_2\in S_Q$.
Note that the definition of the cone is independent of the choice of $x$ for $S_x$, while it plays a role when we consider $T_x$ (this is also the reason why we only allow rational parameters $x$ in the last case).
For the characteristics $a,b\in R(c_1)\cap R(c_2)$, we introduce the holomorphic theta series 
\begin{align}\label{theta_part}
  \Theta_{a,b}^{c_1,c_2}[f] (\tau)=\sum_{n'\in a+\Z^3} \{\sign \bigl(B(n',c_1)\bigr)-\sign \bigl(B(n',c_2)\bigr) \}\ko (n'_1x^2+n'_2x+n'_3)^{k-1}\ko q^{Q(n')} \ko e^{2\pi i B(n',b)}.
\end{align}
Further, we consider $(a,b)$ as the $3\times 2$-matrix with the two column vectors $a$ and $b$ and define the meromorphic function
\begin{align*}
    f_{(a,b)}(\tau):=-\frac{(k-1)!}{(8\pi i)^k}\frac1{(a_1\tau+b_1)^k} \ko e^{-8\pi i a_1b_3}\ko  \sum_{m\in a_2+\Z} q^{m^2}\ko e^{4\pi i b_2m}.
\end{align*}
Also we consider the unary theta function 
\begin{align}\label{eq:unarytheta}
    \theta(\tau):=\sum_{n\in \Z}e^{2\pi in^2\tau}.
\end{align}
The connection with $S_x$ and $T_x$ is then given by:

\begin{lemma}\label{limSxTx}
Let
\begin{align*}
    \widetilde{a}=(a_1,2a_1x+a_2,a_1x^2+a_2x+a_3),\quad\widetilde{b}= (b_1,2b_1x+b_2,b_1x^2+b_2x+b_3)
\end{align*}
and
\begin{align*}
    \widehat{a}=(a_3,2a_3x^{-1}+a_2,a_3x^{-2}+a_2x^{-1}+a_1),\quad\widehat{b}= (b_3,2b_3x^{-1}+b_2,b_3x^{-2}+b_2x^{-1}+b_1).
\end{align*}
(a) Let $c_1=-\frac{1}{4}\Bigl(\begin{smallmatrix}0\\0\\1\end{smallmatrix}\Bigr)$ and $c_2=-\frac{1}{4}\Bigl(\begin{smallmatrix}1\\0\\0\end{smallmatrix}\Bigr)$, both in $S_Q$.
Then we have
\begin{align*}
    S_x(\tau)=\frac 12\lim_{a,b\rightarrow 0}\Bigl(\frac12 \Theta_{a,b}^{c_1,c_2}[f] (\tau)-f_{(\widetilde{a},\widetilde{b})}(\tau)+x^{2k-2}\,f_{(\widehat{a},\widehat{b})}(\tau)\Bigr).
\end{align*}
(b) Let $c_1=-\frac{1}{4}\Bigl(\begin{smallmatrix}1\\-2x\\x^2\end{smallmatrix}\Bigr)$ and $c_2=-\frac{1}{4}\Bigl(\begin{smallmatrix}0\\0\\1\end{smallmatrix}\Bigr)$, both in $S_Q$.
Then we have
\begin{align*}
    T_x(\tau)=\frac 12\lim_{a,b\rightarrow 0}\Bigl(\frac12 \Theta_{a,b}^{c_1,c_2}[f] (\tau)+f_{(\widetilde{a},\widetilde{b})}(\tau)\Bigr)
\end{align*}
for $k\geq 4$.
Let $(a',b')=(a,b)g$, with $g=\left(\begin{smallmatrix}g_{11}&g_{12}\\ g_{21}&g_{22}\end{smallmatrix}\right)$ ($g_{21}$ and $g_{22}$ not both zero).
Then we have
\[T_x(\tau)=\frac12\lim_{b\rightarrow 0} \lim_{a\rightarrow 0}\Bigl(\frac12 \Theta_{a',b'}^{c_1,c_2} [f](\tau)+f_{(\widetilde a', \widetilde b')}(\tau)\Bigr)+\frac{\kappa(x)}{8\pi i} \biggl( \frac12\frac{g_{21}}{g_{21}\tau+g_{22}}\ko \theta(\tau) + \theta'(\tau)\biggr)\]
for $k=2$.
\end{lemma}

To establish the modular behavior of $S_x$ and $T_x$ we first consider the modular behavior of $\Theta_{a,b}^{c_1,c_2}[f]$ and $f_{(a,b)}$.
For this we consider modular substitutions on $\Gamma_0(4)$ and make use of Shimura's definition (see \cite{shimura}) of the automorphic factor $j(\gamma,\tau):=\theta(\gamma \tau)/\theta(\tau)$ for any $\gamma \in \Gamma_0(4)$.
As the modular substitutions alter the characteristics, we introduce the modified characteristics $a',b'$ as the one satisfying $(a',b'):=(a,b)\gamma$.

In the next lemma we show that $f_{(a,b)}$ and the theta function $\Theta_{a,b}^{c_1,c_2}[f]$ defined in \eqref{theta_part} have the same modular transformation behavior on $\Gamma_0(4)$:

\begin{lemma}\label{modthetafab}
We have
\begin{align}\label{align_theta4}
   \Theta_{a,b}^{c_1,c_2}[f](\gamma\tau) =j(\gamma,\tau)^{2k+1}\ko e^{\pi i B(a,b)-\pi i B(a',b')}\ko \Theta_{a',b'}^{c_1,c_2}[f](\tau)
\end{align}
and
\begin{align}\label{align_f}
    f_{(a,b)}(\gamma\tau) =j(\gamma,\tau)^{2k+1}\ko e^{\pi i B(a,b)-\pi i B(a',b')}\ko f_{(a',b')}(\tau)
\end{align}
for all $\gamma\in \Gamma_0(4)$.
\end{lemma}

\begin{proof}[Proof of Theorem \ref{mainprop}]
We note that the maps $(a,b)\mapsto (a',b')=(a,b)\gamma$ and $(a,b)\mapsto (\widetilde a,\widetilde b)$ commute and that $B(\widetilde a',\widetilde b')= B(a',b')$ holds.
Similarly this holds for $\widehat a$ and $\widehat b$.
Hence from Lemma \ref{limSxTx} and Lemma \ref{modthetafab} we directly get
\[ \begin{split}
    S_x(\gamma\tau) &= \frac 12\lim_{a,b\rightarrow 0}\Bigl(\frac12 \Theta_{a,b}^{c_1,c_2}[f] (\gamma\tau)-f_{(\widetilde{a},\widetilde{b})}(\gamma\tau)+x^{2k-2}\,f_{(\widehat{a},\widehat{b})}(\gamma\tau)\Bigr)\\
    &=j(\gamma,\tau)^{2k+1}\ko \frac 12\lim_{a,b\rightarrow 0}e^{\pi i B(a,b)-\pi i B(a',b')}\ko\Bigl(\frac12 \Theta_{a',b'}^{c_1,c_2}[f] (\tau)-f_{(\widetilde{a}',\widetilde{b}')}(\tau)+x^{2k-2}\,f_{(\widehat{a}',\widehat{b}')}(\tau)\Bigr)\\
    &= j(\gamma,\tau)^{2k+1}\ko S_x(\tau)
    \end{split}\]
for all $\gamma\in \Gamma_0(4)$.
In exactly the same way we obtain that for $k\geq 4$
\[ T_x(\gamma\tau)= j(\gamma,\tau)^{2k+1}\ko T_x(\tau)\]
holds for all $\gamma\in \Gamma_0(4)$.
For $k=2$ we use Lemma \ref{limSxTx} twice (first with $g=I$ and then with $g=\gamma=\left(\begin{smallmatrix}\gamma_{11}&\gamma_{12}\\ \gamma_{21}&\gamma_{22}\end{smallmatrix}\right)\in\Gamma_0(4)$) and find
\[\begin{split}
    T_x(\gamma \tau)&=\frac12\lim_{b\rightarrow 0} \lim_{a\rightarrow 0}\Bigl(\frac12 \Theta_{a,b}^{c_1,c_2} [f](\gamma\tau)+f_{(\widetilde a, \widetilde b)}(\gamma\tau)\Bigr)+\frac{\kappa(x)}{8\pi i} \ko \theta'(\gamma\tau)\\
    &=j(\gamma,\tau)^5\ko \frac12\lim_{b\rightarrow 0} \lim_{a\rightarrow 0}e^{\pi i B(a,b)-\pi i B(a',b')}\Bigl(\frac12 \Theta_{a',b'}^{c_1,c_2} [f](\tau)+f_{(\widetilde a', \widetilde b')}(\tau)\Bigr)+\frac{\kappa(x)}{8\pi i} \ko \theta'(\gamma\tau)\\
    &=j(\gamma,\tau)^5\ko T_x(\tau) +\frac{\kappa(x)}{8\pi i} \biggl( \theta'(\gamma\tau) -j(\gamma,\tau)^5 \biggl( \frac12\frac{\gamma_{21}}{\gamma_{21}\tau+\gamma_{22}}\ko \theta(\tau) + \theta'(\tau) \biggr)\biggr).
\end{split}\]
It is known that
\[ j(\gamma,\tau)^2= \Bigl( \frac{-1}{\gamma_{22}}\Bigr) \ko (\gamma_{21}\tau+\gamma_{22}),\]
from which we get
\[ \frac{j'(\gamma,\tau)}{j(\gamma,\tau)}= \frac12\frac{\gamma_{21}}{\gamma_{21}\tau+\gamma_{22}},\]
where $j'$ denotes the derivative of $j$ with respect to $\tau$.
Using $\theta(\gamma\tau)=j(\gamma,\tau)\ko \theta(\tau)$ and $j(\gamma,\tau)^{4}\frac\partial{\partial\tau} (\gamma\tau)=1$ we then obtain
\[ \begin{split}
    T_x(\gamma \tau)&= j(\gamma,\tau)^5\ko T_x(\tau) +\frac{\kappa(x)}{8\pi i} \biggl( \theta'(\gamma\tau) -j(\gamma,\tau)^5 \biggl( \frac{j'(\gamma,\tau)}{j(\gamma,\tau)}\ko \theta(\tau) + \theta'(\tau) \biggr)\biggr)\\
    &= j(\gamma,\tau)^5\ko T_x(\tau) +\frac{\kappa(x)}{8\pi i} \bigl( \theta'(\gamma\tau) -j(\gamma,\tau)^4 \frac{\partial}{\partial\tau} \theta(\gamma\tau)\bigr)= j(\gamma,\tau)^5\ko T_x(\tau),
\end{split}\]
as desired.

Also note that $S_x|\gamma$ and $T_x|\gamma$ are holomorphic at $\infty$ for all $\gamma\in\operatorname{SL}_2(\Z)$: First of all, we know that this holds for $\vartheta$ and $\vartheta'$. Since we can deduce from Theorem 2.5(4) and (5) (analogously to Lemma 3.7 in \cite{RZ}) that this also holds for $\Theta_{a,b}^{c_1,c_2}[f]$ if $a,b\in R(c_1)\cap R(c_2)$ and one can show a similar relation as in the aforementioned lemma for $f_{(a,b)}$, we conclude that $S_x$ and $T_x$ satisfy the desired growth conditions at the cusps of $\Gamma_0(2)$.
\end{proof}

\begin{proof}[Proof of Lemma \ref{limSxTx}]
We consider $\frac12\Theta_{a,b}^{c_1,c_2}[f]$, with the theta series as in \eqref{theta_part}.
To obtain the formula for $S_x$ we set $c_1=-\frac{1}{4}\Bigl(\begin{smallmatrix}0\\0\\1\end{smallmatrix}\Bigr)$ and $c_2=-\frac{1}{4}\Bigl(\begin{smallmatrix}1\\0\\0\end{smallmatrix}\Bigr)$.
Writing $n'=a+\Bigl(\begin{smallmatrix}n_1\\n_2\\n_3\end{smallmatrix}\Bigr)$ with $n_i\in \Z$, this series still converges absolutely if we restrict the summation to the part of the lattice $\Z^3$ for which $n_1\neq 0$ and $n_3\neq 0$ holds and let $a,b\rightarrow 0$.
It doesn't play a role whether we consider the limit of $a$ or of $b$ first, in any case the partial sum asymptotes to
\[\Biggl(\sum_{\substack{n\in \Z^3\\ n_1>0>n_3}}-\sum_{\substack{n\in \Z^3\\ n_1<0<n_3}}\Biggr) (n_1x^2+n_2x+n_3)^{k-1}\ko q^{Q(n)}= 2 \sum_{\substack{n\in \Z^3\\ n_1>0>n_3}} (n_1x^2+n_2x+n_3)^{k-1}\ko q^{Q(n)},\]
substituting $n\mapsto -n$ in the second sum and using the fact that $k-1$ is odd.

We fix $a=\Bigl(\begin{smallmatrix}a_1\\a_2\\a_3\end{smallmatrix}\Bigr)$ and $b=\Bigl(\begin{smallmatrix}b_1\\b_2\\b_3\end{smallmatrix}\Bigr)$, where for simplicity we will only consider the case that $a_1,a_3\in (0,1)$.
The other cases can be dealt with in a similar manner.
In our case we then have $B(c_i,a)> 0$.
We first investigate the series in \eqref{theta_part} on the part, where $n_1=0$ holds:
note that
\[\sign \bigl(B(n',c_1)\bigr)-\sign \bigl(B(n',c_2)\bigr) =\sign (a_1)-\sign (a_3+n_3)\]
equals 2 for strictly negative $n_3$ and vanishes otherwise,
so this partial sum equals
\begin{align}\label{align_n1=0}
\sum_{n_2\in \Z} \sum_{n_3=-\infty}^{-1}  (a_1x^2+a_3+(a_2+n_2)x+n_3)^{k-1}\ko q^{(a_2+n_2)^2-4a_1(a_3+n_3)} \ko e^{-8\pi i ( b_1 (a_3+n_3)+a_1b_3)+4\pi i b_2(a_2+n_2)}.
\end{align}
We set $z_1:=a_1\tau +b_1$ and write \eqref{align_n1=0} as
\begin{multline*}
e^{-8\pi i a_1b_3+8\pi i z_1a_1x^2} \sum_{n_2\in \Z} q^{(a_2+n_2)^2}\ko e^{8\pi i z_1(a_2+n_2)x}\ko  e^{4\pi i b_2(a_2+n_2)}\\
\quad\cdot\sum_{n_3=-\infty}^{-1} (a_1x^2+a_3+(a_2+n_2)x+n_3)^{k-1}\ko e^{-8\pi i z_1(a_1x^2+a_3+(a_2+n_2)x+n_3)}.
\end{multline*}
We make use of Lemma \ref{lemma_bernoulli} with $\alpha=a_1x^2+a_3+(a_2+n_2)x$, $\beta=0$ and $z=-8\pi i z_1$ to rewrite the sum over $n_3$ as
\begin{align*}
 \frac{(-1)^{k-1}\ko (k-1)!}{(-8\pi i z_1)^k}+\sum_{m=0}^\infty \frac{B_{m+k}(a_1x^2+a_3+(a_2+n_2)x)}{m+k}\frac{(-8\pi i z_1)^m}{m!}.   
\end{align*}
First, we consider the second summand, which contains no poles for $a,b\rightarrow 0$: if we let $a,b\rightarrow 0$, this part goes to 
$B_k(n_2x)/k$, as $(-8\pi i z_1)^{m}\rightarrow 0$ whenever $m>0$.
So for the ``regular" part in \eqref{align_n1=0}, we obtain  $\frac{1}{k}\sum_{m\in \Z} B_k(mx)\ko q^{m^2}$ for $a,b\rightarrow 0$.

We define the remaining part of \eqref{align_n1=0}, using that $k$ is even, as
\begin{align*}
f(a,b,x;\tau)&:=-e^{-8\pi i a_1b_3+8\pi i z_1a_1x^2}\ko\frac{(k-1)!}{(8\pi i z_1)^k}  \sum_{m\in a_2+\Z} q^{m^2}\ko e^{8\pi i mxz_1+4\pi ib_2 m}\\
&= -e^{-8\pi i a_1(b_1x^2+b_2x+b_3)} \ko\frac{(k-1)!}{(8\pi i z_1)^k} \sum_{m\in 2a_1x+a_2+\Z} q^{m^2}\ko e^{4\pi i(2b_1x+b_2) m}.
\end{align*}
If we now consider the series in \eqref{theta_part} on the part of the sum, where $n_3=0$ holds, we have
\[-\sum_{n_2\in \Z} \sum_{n_1=-\infty}^{-1}  ((a_1+n_1)x^2+a_3+(a_2+n_2)x)^{k-1}\ko q^{(a_2+n_2)^2-4a_3(a_1+n_1)} \ko e^{-8\pi i ( b_1 a_3+b_3(a_1+n_1))+4\pi ib_2(a_2+n_2)},\]
which equals
\[-x^{2k-2}\sum_{n_2\in \Z} \sum_{n_1=-\infty}^{-1}  \bigl(a_1+n_1+\frac{a_3}{x^2}+\frac{a_2+n_2}{x}\bigr)^{k-1}\ko q^{(a_2+n_2)^2-4a_3(a_1+n_1)} \ko e^{-8\pi i ( b_1 a_3+b_3(a_1+n_1))+4\pi ib_2(a_2+n_2)}.\]
As this series is (up to the factor $-x^{2k-2}$) just \eqref{align_n1=0} for $1/x$ instead of $x$ interchanging $a_1$ and $a_3$ and $b_1$ and $b_3$, respectively, we obtain the function $-x^{2k-2}\ko f\Bigl(\Bigl(\begin{smallmatrix}a_3\\a_2\\a_1\end{smallmatrix}\Bigr),\Bigl(\begin{smallmatrix}b_3\\b_2\\b_1\end{smallmatrix}\Bigr),\frac{1}{x};\tau\Bigr)$ and as ``regular" part 
\[ -x^{2k-2} \frac{1}{k}\sum_{m=-\infty}^\infty B_k\bigl(\frac mx\bigr)\,q^{m^2}.\]
Combining these results for the different parts of the lattice $\Z^3$ we then get
\begin{align*}
    \lim_{a,b\rightarrow 0}&\Bigl(\frac12 \Theta_{a,b}^{c_1,c_2} (\tau)-f(a,b,x;\tau)+x^{2k-2}f\Bigl(\Bigl(\begin{smallmatrix}a_3\\a_2\\a_1\end{smallmatrix}\Bigr),\Bigl(\begin{smallmatrix}b_3\\b_2\\b_1\end{smallmatrix}\Bigr),\frac{1}{x};\tau\Bigr)\Bigr)\\
    &= 2 \sum_{\substack{n\in \Z^3\\ n_1>0>n_3}} (n_1x^2+n_2x+n_3)^{k-1}\ko q^{Q(n)} +\frac{1}{k}\sum_{m=-\infty}^\infty B_k(mx)\ko q^{m^2} -x^{2k-2}\frac{1}{k}\sum_{m=-\infty}^\infty B_k\bigl(\frac{m}{x}\bigr)\ko q^{m^2}\\
    &= 2 S_x(\tau).
\end{align*}
In terms of $\widetilde a$ and $\widetilde b$ we have $f(a,b,x;\tau)=f_{(\widetilde{a},\widetilde{b})}(\tau)$ and similarly 
$f\Bigl(\Bigl(\begin{smallmatrix}a_3\\a_2\\a_1\end{smallmatrix}\Bigr),\Bigl(\begin{smallmatrix}b_3\\b_2\\b_1\end{smallmatrix}\Bigr),\frac{1}{x};\tau\Bigr)=f_{(\widehat{a},\widehat{b})}(\tau)$.
This gives the desired result for $S_x$.

To get the formula for $T_x$ we now set $c_1=-\frac{1}{4}\Bigl(\begin{smallmatrix}1\\-2x\\x^2\end{smallmatrix}\Bigr)$ and $c_2=-\frac{1}{4}\Bigl(\begin{smallmatrix}0\\0\\1\end{smallmatrix}\Bigr)$.
We then have
\[ \begin{split}
\frac12 \Theta_{a,b}^{c_1,c_2}[f](\tau)= \frac 12\sum_{n\in\Z^3} \bigl\{ \sign(n_1x^2&+n_2x+n_3+\widetilde{a}_3) -\sign(n_1+a_1)\bigr\}(n_1x^2+n_2x+n_3+\widetilde{a}_3)^{k-1} \\
&\cdot q^{(n_2+a_2)^2-4(n_1+a_1)(n_3+a_3)}\ko e^{4\pi i((n_2+a_2)b_2-2(n_1+a_1)b_3-2(n_3+a_3)b_1)}.
\end{split}\]
Similar to the previous case we assume $\widetilde{a}_3,a_1\in(0,1)$ for simplicity and first consider the part of $\Z^3$ where $n_1x^2+n_2x+n_3\neq 0$ and $n_1\neq 0$ hold.
Letting $a,b\rightarrow 0$ we then obtain
\[\begin{split}\Biggl(\sum_{\substack{n\in \Z^3\\ n_1x^2+n_2x+n_3>0>n_1}}&-\sum_{\substack{n\in \Z^3\\ n_1x^2+n_2x+n_3<0<n_1}}\Biggr) (n_1x^2+n_2x+n_3)^{k-1}\ko q^{Q(n)}\\
&= 2 \sum_{\substack{n\in \Z^3\\ n_1x^2+n_2x+n_3>0>n_1}} (n_1x^2+n_2x+n_3)^{k-1}\ko q^{Q(n)}.
\end{split}\]
For the part of $\Z^3$ where $n_1=0$ holds we get
\[\begin{split}
-&\sum_{\substack{n_2,n_3\in\Z\\ n_2x+n_3<0}} (n_2x+n_3+\widetilde{a}_3)^{k-1}\ko  q^{(n_2+a_2)^2-4a_1(n_3+a_3)}\ko e^{4\pi i((n_2+a_2)b_2-2a_1b_3-2(n_3+a_3)b_1)}\\
&= -e^{-8\pi ia_1\widetilde{b}_3} \sum_{\substack{n_2,n_3\in\Z\\ n_2x+n_3<0}} (n_2x+n_3+\widetilde{a}_3)^{k-1}\ko q^{(n_2+\widetilde{a}_2)^2}\ko e^{4\pi i(n_2+\widetilde{a}_2)\widetilde{b}_2}\ko e^{-8\pi iz_1(n_2x+n_3+\widetilde{a}_3)},
\end{split}\]
where again $z_1=a_1\tau+b_1$.
Using Lemma \ref{lemma_bernoulli} with $\alpha=\widetilde{a}_3$, $\beta=n_2x$ and $z=-8\pi i z_1$, we get
\[-\sum_{\substack{n_3\in\Z\\ n_2x+n_3<0}} (n_2x+n_3+\widetilde{a}_3)^{k-1}\ko e^{-8\pi iz_1(n_2x+n_3+\widetilde{a}_3)}=\frac{(k-1)!}{(8\pi iz_1)^k} - \sum_{m=0}^\infty \frac{B_{m+k}(\widetilde{a}_3+n_2x-\lfloor n_2x\rfloor)}{m+k}\frac{(-8\pi iz_1)^m}{m!}\]
so for the ``regular'' part with $a,b\rightarrow 0$ we have
\[-\sum_{n_2\in\Z} q^{n_2^2}\frac{B_k(n_2x-\lfloor n_2x \rfloor)}{k}=-\frac1k \sum_{m\in\Z} \overline{B}_k(mx)\ko q^{m^2}. \]
The remaining part is
\[ e^{-8\pi ia_1\widetilde{b}_3}\frac{(k-1)!}{(8\pi iz_1)^k} \sum_{n_2\in\Z} q^{(n_2+\widetilde{a}_2)^2}\ko e^{4\pi i(n_2+\widetilde{a}_2)\widetilde{b}_2}=-f(a,b,x;\tau)=-f_{(\widetilde a,\widetilde b)}(\tau).\]
Finally, we consider the part of $\Z^3$ where $n_1x^2+n_2x+n_3=0$ holds and get
\[ \widetilde{a}_3^{k-1} \sum_{\substack{n_1,n_2,n_3\in\Z\\ n_1x^2+n_2x+n_3=0>n_1}} q^{Q(n+a)}\ko e^{2\pi iB(n+a,b)}=\widetilde{a}_3^{k-1} q^{Q(a)}\ko e^{2\pi iB(a,b)}\sum_{\substack{n_1,n_2,n_3\in\Z\\ n_1x^2+n_2x+n_3=0>n_1}} q^{Q(n)}\ko e^{2\pi iB(n,z)},\]
with $z:=a\tau +b$.
We write $x=\frac rs$ with $s\in\N$, $\gcd(r,s)=1$ and consider the substitution
\[ (m_1,m_2,m_3)= (\frac1s n_1, 2\frac rs n_1 +n_2, \frac{r^2}s n_1+rn_2+sn_3).\]
It gives a bijection between the sets $\{ n\in\Z^3 \mid n_1 \equiv 0 \smod s\}$ and $\{ m\in\Z^3 \mid r^2 m_1-r m_2+m_3\equiv 0 \smod s\}$.
We note that if $n_1x^2+n_2x+n_3=0$ holds, then we have $n_1\equiv 0\smod s$.
Further, we have $Q(n)=Q(m)$ and $B(n,z)= 2\widetilde z_2 m_2 - 4\frac1s \widetilde z_1m_3 -4 s \widetilde z_3 m_1$, with $\widetilde z=(z_1, 2z_1x +z_2, z_1x^2+z_2x+z_3)$,
so with this substitution we obtain
\begin{equation}\label{partk2}
\begin{split}
\widetilde{a}_3^{k-1} & q^{Q(a)}\ko e^{2\pi iB(a,b)}\sum_{\substack{m_1,m_2,m_3\in\Z\\ m_3=0>m_1,\ r^2 m_1-r m_2+m_3\equiv 0 \smod s}} q^{Q(m)}\ko e^{4\pi i(\widetilde z_2 m_2 - 2\frac1s \widetilde z_1m_3 -2 s \widetilde z_3 m_1)}\\
&=\widetilde{a}_3^{k-1} q^{Q(a)}\ko e^{2\pi iB(a,b)}\sum_{m_2\in\Z}q^{m_2^2}\ko e^{4\pi i\widetilde z_2 m_2 }\sum_{\substack{m_1\in\Z\\ m_1<0,\ m_1\equiv r^* m_2\smod s}} e^{-8\pi is \widetilde z_3 m_1}\\
&=\frac{\widetilde{a}_3^{k-1} q^{Q(a)}\ko e^{2\pi iB(a,b)}}{e^{-8\pi i \widetilde z_3 s^2}-1} \sum_{m_2\in\Z}q^{m_2^2}\ko e^{4\pi i\widetilde z_2 m_2-8\pi i \widetilde z_3 s^2(\frac{r^*}s m_2 -\lfloor\frac{r^*}s m_2\rfloor)},
\end{split}
\end{equation}
where $r^*$ is the multiplicative inverse of $r$ modulo $s$ and we used the geometric series in the last step.
As in the proof of Lemma \ref{prop_eisenstein} this part vanishes for $k\geq 4$ if we let $a,b\rightarrow 0$.
Combining the results for the different parts we thus find for $k\geq 4$:
\begin{align*}
    \lim_{a,b\rightarrow 0}&\Bigl(\frac12 \Theta_{a,b}^{c_1,c_2} [f](\tau)+f_{(\widetilde a, \widetilde b)}(\tau)\Bigr)\\
    &= 2 \sum_{\substack{n\in \Z^3\\ n_1x^2+n_2x+n_3>0>n_1}} (n_1x^2+n_2x+n_3)^{k-1}\ko q^{Q(n)} - \frac{1}{k}\sum_{m=-\infty}^\infty \overline B_k(mx)\ko q^{m^2}= 2 T_x(\tau),
\end{align*}
which gives the desired result.
For $k=2$ we have to be careful about the order of the limits.
We denote the last expression in \eqref{partk2} by $g_{a,b}(\tau)$ and consider 
\[\lim_{b\rightarrow 0} \lim_{a\rightarrow 0} g_{a',b'}(\tau),\]
where $(a',b')=(a,b)g$, with $g=\left(\begin{smallmatrix}g_{11}&g_{12}\\ g_{21}&g_{22}\end{smallmatrix}\right)$ ($g_{21}$ and $g_{22}$ not both zero).
We find
\[ \begin{split}
\lim_{b\rightarrow 0} \lim_{a\rightarrow 0} g_{a',b'}(\tau) &= \sum_{m_2\in\Z} q^{m_2^2} \cdot \lim_{b\rightarrow 0} \lim_{a\rightarrow 0}\frac{\widetilde a_3'}{e^{-8\pi i \widetilde z'_3 s^2}-1}  \\
&= \sum_{m_2\in\Z} q^{m_2^2} \cdot \lim_{b\rightarrow 0} \lim_{a\rightarrow 0}\frac{g_{11}\widetilde a_3+g_{21}\widetilde b_3}{e^{-8\pi i ((g_{11}\widetilde a_3 +g_{21} \widetilde b_3)\tau + g_{12}\widetilde a_3 +g_{22} \widetilde b_3) s^2}-1}\\
&= \sum_{m_2\in\Z} q^{m_2^2} \cdot \lim_{b\rightarrow 0}\frac{g_{21}\widetilde b_3}{e^{-8\pi i (g_{21} \tau + g_{22} ) \widetilde b_3s^2}-1} \\
&= -\frac1{8\pi i} \frac1{s^2}\frac{g_{21}}{g_{21}\tau+g_{22}}\sum_{m_2\in\Z} q^{m_2^2} = -\frac{\kappa(x)}{8\pi i} \frac{g_{21}}{g_{21}\tau+g_{22}}\ko \theta(\tau).
\end{split}\]
Combining the results for the different parts we thus find for $k=2$:
\begin{align*}
     \lim_{b\rightarrow 0} &\lim_{a\rightarrow 0}\Bigl(\frac12 \Theta_{a',b'}^{c_1,c_2} [f](\tau)+f_{(\widetilde a', \widetilde b')}(\tau)\Bigr)\\
    &= 2 \sum_{\substack{n\in \Z^3\\ n_1x^2+n_2x+n_3>0>n_1}} (n_1x^2+n_2x+n_3)\ko q^{Q(n)} - \frac{1}{2}\sum_{m=-\infty}^\infty \overline B_2(mx)\ko q^{m^2}-\frac{\kappa(x)}{8\pi i} \frac{g_{21}}{g_{21}\tau+g_{22}}\ko \theta(\tau)\\
    &=2 T_x(\tau) -\frac{\kappa(x)}{8\pi i} \frac{g_{21}}{g_{21}\tau+g_{22}}\ko \theta(\tau)-\kappa(x) \sum_{m=1}^\infty m^2 q^{m^2}=2 T_x(\tau) -\frac{\kappa(x)}{4\pi i} \biggl( \frac12\frac{g_{21}}{g_{21}\tau+g_{22}}\ko \theta(\tau) + \theta'(\tau)\biggr),
\end{align*}
which gives the desired result for this case.
\end{proof}

\begin{proof}[Proof of Lemma \ref{modthetafab}]
We can easily see that \eqref{align_theta4} and \eqref{align_f} are compatible with matrix multiplication.
Hence we only have to show them for the generators $-I,T$ and $\gamma'=\bigl(\begin{smallmatrix}1&0\\4&1\end{smallmatrix}\bigr)$ of $\Gamma_0(4)$. Note that $j(-I,\tau)=j(T,\tau)=1$ and $j(\gamma',\tau)^2=4\tau+1$. We also note that $(a,b)\ko (-I)=-(a,b)$, $(a,b)\ko T=(a,a+b)$ and $(a,b)\ko\gamma'=(a+4b,b)$.

For the sake of better readability we write $\Theta_{(a,b)}$ instead of $\Theta_{a,b}^{c_1,c_2}[f]$. 
By Theorem \ref{modularity}(3), $\Theta_{-(a,b)}=\Theta_{(a,b)}$ holds since the degree of $f$ is odd. Thus we have
\[\Theta_{(a,b)}(-I\tau)=\Theta_{(a,b)}(\tau)=\Theta_{-(a,b)}(\tau),\]
which is \eqref{align_theta4} for $\gamma=-I$.
From Theorem \ref{modularity}(4) we obtain 
\[\Theta_{(a,b)}(\tau+1)=e^{-2\pi i Q(a)}\ko\Theta_{(a,a+b)}(\tau),\]
which is \eqref{align_theta4} for $\gamma=T$.
As a result, we also have
\begin{align}\label{align_TN}
    \Theta_{(a,b)}(\tau-4)=e^{8\pi i Q(a)}\ko\Theta_{(a,-4a+b)}(\tau).
\end{align}
By Theorem \ref{modularity}(5) we have
\begin{align*}
    \Theta_{(a,b)}\Bigl(-\frac{1}\tau\Bigr)=(-i\tau)^{k+1/2}\ko \frac{i^k}{4\sqrt{2}}\ko e^{2\pi i B(a,b)}\ko \sum_{p\in A^{-1}\Z^3/\Z^3} \Theta_{(b+p,-a)}(\tau).
\end{align*}
To keep track of the choice of the square root, we write $\sqrt{-i\tau}$ as
\[ \sqrt{2}\, \frac{\theta(-1/\tau)}{\theta(\tau/4)},\]
where $\theta$ is defined as in \eqref{eq:unarytheta}.
Hence we have
\begin{equation}\label{eq:trans} \frac{\Theta_{(a,b)}}{\theta^{2k+1}}\Bigl(-\frac{1}\tau\Bigr)=2^{k-2} i^k\ko e^{2\pi i B(a,b)}\ko \frac1{\theta(\tau/4)^{2k+1}}\sum_{p\in A^{-1}\Z^3/\Z^3} \Theta_{(b+p,-a)}(\tau).
\end{equation}
Using \eqref{align_TN} we obtain
\[ \Theta_{(b+p,-a)}(\tau -4)=e^{8\pi iQ(b+p)}\ko \Theta_{(b+p,-4(b+p)-a)}(\tau).\]
Next we use Theorem \ref{modularity}(2) with $\mu=-4p$.
Note that $4A^{-1}$ is even, so we have $\mu\in\Z^3$ and $4Q(p)\in\Z$.
Hence we get
\begin{align*}
    \Theta_{(b+p,-a)}(\tau -4)&=e^{8\pi i Q(b+p)-8\pi i B(b+p,p)}\ko\Theta_{(b+p,-a-4b)}(\tau)\\
    &=e^{8\pi iQ(b)-8\pi i Q(p)}\ko \Theta_{(b+p,-a-4b)}(\tau)\\
    &=e^{8\pi iQ(b)}\ko \Theta_{(b+p,-a-4b)}(\tau).
\end{align*}
We replace $\tau$ by $\tau-4$ in \eqref{eq:trans}, use the identity above and then \eqref{eq:trans} again with $a$ replaced by $a+4b$ to obtain
\[\frac{\Theta_{(a,b)}}{\theta^{2k+1}}\Bigl(-\frac{1}{\tau-4}\Bigr) = e^{-8\pi iQ(b)}\frac{\Theta_{(a+4b,b)}}{\theta^{2k+1}}\Bigl(-\frac{1}\tau\Bigr).\]
Replacing $\tau$ by $-1/\tau$ then gives \eqref{align_theta4} for $\gamma=\gamma'=\bigl(\begin{smallmatrix}1&0\\4&1\end{smallmatrix}\bigr)$.

The function $f_{(a,b)}$ is the product of a function $g_{a_1,b_1}$ that transforms like a modular form of weight $k$ and a function $\theta_{a_2,b_2}$ that transforms like a modular form of weight $1/2$, as we show in the following.
 We define
\[g_{a_1,b_1}(\tau):=\frac{1}{z_1^k}=\frac{1}{(a_1\tau +b_1)^k}\quad\text{and}\quad \theta_{a_2,b_2}(\tau):=\sum_{m\in a_2+\Z}q^{m^2}\ko e^{4\pi i b_2m}.\]
We can easily check that $\theta_{-a_2,-b_2}=\theta_{a_2,b_2}$ holds, which implies $f_{-(a,b)}=f_{(a,b)}$. Thus we have
\[f_{(a,b)}(-I\tau)=f_{(a,b)}(\tau)=f_{-(a,b)}(\tau),\]
which is \eqref{align_f} for $\gamma=-I$.

By straightforward calculations the modular transformation behavior of $g_{a_1,b_1}$ for $T$ and $\gamma'$ follows:
\begin{align*}
    g_{a_1,b_1}(\tau+1)&=\frac{1}{(a_1\tau +a_1+b_1)^k}=g_{a_1,a_1+b_1}(\tau),\\
    g_{a_1,b_1}\Bigl(\frac{\tau}{4\tau+1}\Bigr)&=\frac{(4\tau+1)^k}{\bigl((a_1+4b_1)\tau +b_1\bigr)^k}=j(\gamma',\tau)^{2k}\ko g_{a_1+4b_1,b_1}(\tau),
\end{align*}
where we use in the last equation $j(\gamma',\tau)^{2}=4\tau +1$.

For $\theta_{a_2,b_2}$, we can directly compute for any $N\in \Z$
\begin{align}\label{eq:trans1}
\theta_{a_2,b_2}(\tau+N)=\sum_{m\in a_2+\Z}q^{m^2}\ko e^{4\pi i b_2m+2\pi i m^2N}=e^{-2\pi i Na_2^2}\ko \theta_{a_2,b_2+Na_2}(\tau),
\end{align}
as $m^2N\equiv 2Na_2m-Na_2^2\ko(\smod{\Z})$.

Since we know that the classical Jacobi theta function $\oldtheta_{a,b}(\tau):=\sum_{m\in a+\Z} e^{\pi i m^2\tau +2\pi i bm}$ satisfies
\[\oldtheta_{a,b}\Bigl(-\frac{1}{\tau}\Bigr)=e^{2\pi i ab}\ko\sqrt{-i\tau}\ko \oldtheta_{b,-a}(\tau),\]
we derive the functional equation
\[\theta_{a_2,b_2}\Bigl(-\frac{1}{\tau}\Bigr)=e^{4\pi i a_2b_2}\ko\sqrt{\frac{-i\tau}{2}}\ko \theta_{-2b_2,\frac{a_2}{2}}\Bigl(\frac{\tau}{4}\Bigr)\]
and
\begin{align}\label{eq:trans2}
\frac{\theta_{a_2,b_2}}{\theta}\Bigl(-\frac{1}{\tau}\Bigr)=e^{4\pi i a_2b_2}\ko \frac{\theta_{-2b_2,\frac{a_2}{2}}}{\theta}\Bigl(\frac{\tau}{4}\Bigr).
\end{align}
Replacing $\tau$ by $\tau-4$ in this equation and applying \eqref{eq:trans1} for $N=-1$, we obtain 
\[\frac{\theta_{a_2,b_2}}{\theta}\Bigl(-\frac{1}{\tau-4}\Bigr)=e^{4\pi i a_2b_2+8\pi i b_2^2}\ko \frac{\theta_{-2b_2,\frac{a_2}{2}+2b_2}}{\theta}\Bigl(\frac{\tau}{4}\Bigr).\]
We use \eqref{eq:trans2} again with $a_2$ replaced by $a_2+4b_2$ and get
\[\frac{\theta_{a_2,b_2}}{\theta}\Bigl(-\frac{1}{\tau-4}\Bigr)=e^{-8\pi ib_2^2}\ko\frac{\theta_{a_2+4b_2,b_2}}{\theta}\Bigl(-\frac{1}{\tau}\Bigr).\]
Replacing $\tau$ by $-1/\tau$ then gives
\[ \theta_{a_2,b_2}(\gamma'\tau)=e^{-8\pi i b_2^2}\ko j(\gamma',\tau)\ko \theta_{a_2+4b_2,b_2}(\tau).\]
Now it just remains to be shown that for $(a',b')=(a,b)\ko T= (a,a+b)$ we have
\[\pi i B(a,b)-\pi i B(a',b')-8\pi i a_1'b_3'=-8\pi i a_1b_3-2\pi i a_2^2\]
and for $(a',b')=(a,b)\ko \gamma'= (a+4b,b)$ we have
\[\pi i B(a,b)-\pi i B(a',b')-8\pi i a_1'b_3'=-8\pi i a_1b_3-8\pi i b_2^2.\]
These identities follow by straightforward calculations and thus \eqref{align_f} holds for $-I,T$ and $\gamma'$ and therefore on the whole group $\Gamma_0(4)$.
\end{proof}

\subsection{Hurwitz class numbers}
Chen and Garvan \cite{CG} investigate generating functions $\mathcal{H}_{a,b}$ for Hurwitz class numbers $H(an+b)$ with $a\steilt 24$ and $(a,b)=1$ and show, inter alia, that $\mathcal{H}_{8,7}(q)\equiv \frac{A(-q)}{-q} (\smod 4)$, where $A$ is a certain second order mock theta function. They use the following identity shown by Humbert \cite{humbert}:
\[\mathcal{H}_{8,7}(q)=\frac{1}{q(q)^3_{\infty}} \sum_{m=0}^\infty \frac{(-1)^{m+1}m^2q^{m(m+1)/2}}{1+q^m}\quad\text{with}\quad (q)_\infty :=\prod_{m=1}^\infty (1-q^m),\]
which suggests that we can interpret the sum that defines $\mathcal{H}_{8,7}$ as the holomorphic part of a theta series of exactly the type that we investigated throughout this paper. Indeed, taking $A=\bigl(\begin{smallmatrix}1&1\\1&0\end{smallmatrix}\bigr)$, $f(n)=n_1^2$ (which is a spherical polynomial with regard to $Q$), $a=\frac{1}{2}(\begin{smallmatrix}0\\1\end{smallmatrix})$, $b=\frac{1}{2}(\begin{smallmatrix}1\\0\end{smallmatrix})$, and $c_1=(\begin{smallmatrix}0\\1\end{smallmatrix})\in S_Q$, $c_2=\sqrt{2}(\begin{smallmatrix}-1\\1\end{smallmatrix})\in C_Q$, we obtain the theta series
\begin{align}\label{theta_hurwitz}
  \modtheta_{a,b}^{c_1,c_2}[f](\tau) = y^{-1} \sum_{n\in \frac{1}{2}(\begin{smallmatrix}0\\1\end{smallmatrix})+\Z^2} \bigl\{\sign (n_1) f(ny^{1/2}) - p^{c_2}[f](n y^{1/2})\bigr\}\ko q^{\frac12 n_1^2+n_1n_2}\ko e^{\pi i(n_1+n_2)}.  
\end{align}
Note that $a\notin R(c_1)$, but the series is still well-defined, as the summand for $n_1=0$ vanishes.
In the following theorem we give the connection between $\mathcal{H}_{8,7}$ and $\modtheta_{a,b}^{c_1,c_2}[f]$, and give their modular behavior in terms of the classical Jacobi theta function
\[\oldtheta_2(\tau)=\sum_{n\in \Z} e^{\pi i(n+\frac12)^2 \tau}.\]

\begin{theorem}
The function $\tau\mapsto q^{7/8}\,\mathcal{H}_{8,7}(q)$ is a mock theta function of weight $3/2$ with shadow proportional to $\oldtheta_2$: it is the holomorphic part of the harmonic Maass form
\[\mathcal{F}:=\frac i4\frac{\modtheta_{a,b}^{c_1,c_2}[f]}{\eta^3}\quad\text{with }\modtheta_{a,b}^{c_1,c_2}[f]\text{ as in \eqref{theta_hurwitz}}\text{ and }\eta(\tau):=q^{1/24}(q)_\infty.\]
In particular, we have
\[\mathcal{F}(\tau)= q^{7/8}\,\mathcal{H}_{8,7}(q)+\frac{1}{4\pi i}\int_{-\overline{\tau}}^{i\infty} \frac{\oldtheta_2(z)}{(-i(z+\tau))^{3/2}}\,dz,\]
where $\mathcal{F}$ is a harmonic Maass form of weight $3/2$ on the congruence subgroup $\Gamma_0(2)$.
The transformation properties of $\mathcal{F}$ are such that $\oldtheta_2 \mathcal F$ transforms as a modular form of weight 2 on $\Gamma_0(2)$ (without a character).
\end{theorem}

\begin{proof}
First we consider the modular transformation behavior of $\mathcal{F}$:
from Theorem \ref{modularity} we get
\[\modtheta_{a,b}^{c_1,c_2}[f](\tau+1)=\modtheta_{a,b}^{c_1,c_2}[f](\tau)\qquad\text{and}\qquad \modtheta_{a,b}^{c_1,c_2}[f]\Bigl(-\frac{1}{\tau}\Bigr)=i\tau^3\,\modtheta_{b,-a}^{c_1,c_2}[f](\tau).\]
Using the last equation (twice) we also find
\[ \modtheta_{a,b}^{c_1,c_2}[f]\Bigl(\frac{\tau}{2\tau+1}\Bigr)= -i(2\tau+1)^3\, \modtheta_{a,b}^{c_1,c_2}[f](\tau).\]
Using the well-known transformation behavior of $\eta$ and $\oldtheta_2$ we get
\[ \frac{\oldtheta_2}{\eta^3} (\tau+1)= \frac{\oldtheta_2}{\eta^3} (\tau)\qquad\text{and}\qquad \frac{\oldtheta_2}{\eta^3} \Bigl(\frac{\tau}{2\tau+1}\Bigr)= \frac i{2\tau+1}\ko \frac{\oldtheta_2}{\eta^3} (\tau)\]
and so $\oldtheta_2 \mathcal{F} = (i/4)\ko \oldtheta_2/\eta^3\cdot \modtheta_{a,b}^{c_1,c_2}[f]$ satisfies
\[ (\oldtheta_2 \mathcal{F})(\tau+1)= (\oldtheta_2 \mathcal{F})(\tau)\qquad\text{and}\qquad (\oldtheta_2 \mathcal{F})\Bigl(\frac{\tau}{2\tau+1}\Bigr)= (2\tau+1)^2 (\oldtheta_2 \mathcal{F})(\tau).\]
Since $\Gamma_0(2)$ is generated by $-I$, $T$ and $\left(\begin{smallmatrix}1&0\\2&1\end{smallmatrix}\right)$ we thus get that
$\oldtheta_2 \mathcal{F}$ transforms as a modular form of weight 2 on $\Gamma_0(2)$ (without character).
Hence
\[\mathcal{F}\Bigl(\frac{a\tau+b}{c\tau+d}\Bigr)=\zeta(\gamma)\,(c\tau+d)^{3/2}\,\mathcal{F}(\tau)\quad\text{for }\gamma=\bigl(\begin{smallmatrix}a&b\\c&d\end{smallmatrix}\bigr)\in \Gamma_0(2),\]
where $\zeta(\gamma)=\sqrt{c\tau+d}\ko \oldtheta_2(\tau)/\oldtheta_2(\gamma\tau)$ is an eighth root of unity.

By the definition of $p^{c}[f]$ and using $E(z)=\sign(z)(1-\beta(z^2))$ with $\beta(x)=\int_x^\infty u^{-1/2} e^{-\pi u} du$, we have 
\begin{multline*}
  y^{-1} \{\sign (v_1) f(vy^{1/2}) - p^{c_2}[f](vy^{1/2})\}\\
  =\{\sign (v_1)+\sign (v_2)\} v_1^2+\Big\{-\sign(v_2)\beta(2yv_2^2)v_1^2-\frac{1}{\pi\sqrt{2y}}v_1E'\bigl(v_2\sqrt{2y}\bigr)+\frac{1}{8\pi^2y}E''\bigl(v_2\sqrt{2y}\bigr)\Big\}.
\end{multline*}
This already shows how we can decompose the theta function into its holomorphic and non-holomorphic part.
For the holomorphic part we immediately derive:
\begin{align*}
i\sum_{m,n\in\Z} & \{\sign (m)+\sign (n+1/2)\}\, (-1)^{m+n}\ko m^2\, q^{\frac12 m^2+m(n+\frac12)}\\
&= i \sum_{m\in\Z\smallsetminus\{0\}} (-1)^m\ko m^2\ko q^{m(m+1)/2} \sum_{n\in\Z}\{\sign (m)+\sign (n+1/2)\}\, (-q^m)^n\\
&=2i \sum_{m\in\Z\smallsetminus\{0\}}\frac{(-1)^m\, m^2\, q^{m(m+1)/2}}{1+q^m}
    =-4i\sum_{m=1}^\infty \frac{(-1)^{m+1}\, m^2\, q^{m(m+1)/2}}{1+q^m}
\end{align*}
So $q^{7/8}\,\mathcal{H}_{8,7}(q)$ is the holomorphic part of $\mathcal{F}$ since $q(q)^3_\infty=q^{7/8}\eta^3(\tau)$.

For the non-holomorphic part of $\modtheta_{a,b}^{c_1,c_2}[f]$ we have
\begin{equation}\label{align_nh}
\begin{split}
   i\sum_{m\in \Z}\sum_{n\in \Z} \Big\{-\sign\bigl(n+\frac12\bigr)\beta\bigl(2y(n+\frac12)^2\bigr)m^2&-\frac{1}{\pi\sqrt{2y}}E'\bigl((n+\frac12)\sqrt{2y}\bigr)m\\
   &+\frac{1}{8\pi^2y}E''\bigl((n+\frac12)\sqrt{2y}\bigr)\Big\}
    \, (-1)^{m+n} q^{m(m+1)/2 +mn}.
    \end{split}
\end{equation}
When we substitute $m-n$ for $m$, \eqref{align_nh} can be written as the sum of the three series
\begin{align}
    -i&\sum_{m\in \Z}\sum_{n\in \Z} \sign\bigl(n+\frac12\bigr)\beta\bigl(2y(n+\frac12)^2\bigr)(m-n)^2\,(-1)^{m} q^{m(m+1)/2 -n(n+1)/2},\label{nh1}\\ 
    \frac{-i}{\pi \sqrt{2y}}&\sum_{m\in \Z}\sum_{n\in \Z}E'\bigl((n+\frac12)\sqrt{2y}\bigr)\,(m-n)\,(-1)^{m} q^{m(m+1)/2-n(n+1)/2},\label{nh2}\\ 
    \frac{i}{8\pi^2y}&\sum_{m\in \Z}\sum_{n\in \Z} E''\bigl((n+\frac12)\sqrt{2y}\bigr)\,(-1)^{m} q^{m(m+1)/2-n(n+1)/2}.\label{nh3}
\end{align}
From the identities 
\[\sum_{m\in \Z}(-1)^m q^{m(m+1)/2}=0,\quad \sum_{m\in \Z}(-1)^m m q^{m(m+1)/2}=(q)_{\infty}^3\quad\text{and}\quad \sum_{m\in \Z}(-1)^m m^2q^{m(m+1)/2}=-(q)_{\infty}^3,\]
we get that \eqref{nh1} equals 
\[i(q)_{\infty}^3 \sum_{n\in \Z}|2n+1|\,\beta\bigl(2(n+\frac12)^2y\bigr) q^{-n(n+1)/2}.\]
For \eqref{nh2} we also use $E'(z)=2e^{-\pi z^2}$ and obtain
\[ \frac{-i}{\pi} (q)_{\infty}^3 \sqrt{\frac{2}{y}} \sum_{n\in \Z} e^{-2\pi (n+\frac12)^2y}\, q^{-n(n+1)/2}.\]
Furthermore, we see that \eqref{nh3} vanishes.
As $\beta$ is related to the incomplete gamma function, we define $\beta(\alpha;x):=\int_x^{\infty} u^{\alpha -1}e^{-\pi u}\,du$, i.\,e. $\beta(x)=\beta(\frac12;x)$. By partial integration we then see that
\[\beta\bigl(\frac12;x\bigr)=\frac{1}{\pi\sqrt{x}}e^{-\pi x}-\frac{1}{2\pi}\beta\bigl(-\frac12;x\bigr).\]
So \eqref{align_nh} can be simplified to
\begin{align*}
    \frac{-i}{2\pi} (q)_{\infty}^3 \sum_{n\in \Z} |2n+1|\beta\bigl(-\frac12;2(n+\frac12)^2y\bigr)\,q^{-n(n+1)/2}.
\end{align*}
By the substitution of $u$ by $(n+\frac12)^2 u$, we have
\[|2n+1| \beta\bigl(-\frac12;2(n+\frac12)^2y\bigr)=|2n+1|\int_{2(n+\frac12)^2y}^\infty u^{-3/2}e^{-\pi u}\,du=2 \int_{2y}^\infty u^{-3/2}e^{-\pi(n+\frac12)^2 u}\,du.\]
Now we can write the non-holomorphic part of the theta function as
\begin{align*}
    \frac{-i}{\pi} q^{1/8} (q)_{\infty}^3  \int_{2y}^\infty \frac{\sum_{n\in \Z} e^{-\pi(n+\frac12)^2 u}q^{-(n+\frac12)^2/2}}{u^{3/2}}\,du=-\frac{1}{\pi}\eta^3(\tau)\int_{-\overline{\tau}}^{i\infty} \frac{\sum_{n\in \Z} e^{\pi i(n+\frac12)^2 z}}{(-i(z+\tau))^{3/2}}\,dz,
\end{align*}
where in the last step we have substituted $u=-i(z+\tau)$.
Combining the holomorphic and non-holomorphic parts we obtain the desired result:
\[\mathcal{F}(\tau)=q^{7/8}\,\mathcal{H}_{8,7}(q)+\frac{1}{4\pi i}\int_{-\overline{\tau}}^{i\infty} \frac{\oldtheta_2(z)}{(-i(z+\tau))^{3/2}}\,dz\]
If we apply the usual $\xi$-operator, given by
\[ \xi_k(f)(\tau)= 2i y^k\ko \overline{\frac{\partial f}{\partial \overline{\tau}}},\]
in weight $k=3/2$ we directly get
\[ \xi_{3/2}(\mathcal{F})(\tau)= -\frac1{4\pi \sqrt{2}}\ko \oldtheta_2(\tau),\]
which is a holomorphic function.
Hence it's annihilated by $\xi_{1/2}$:
\[ \xi_{1/2}(\xi_{3/2}(\mathcal{F})) =0\]
Since the weight $k$ hyperbolic Laplacian
\[ \Delta_k=-y^2\Bigl( \frac{\partial^2}{\partial x^2}+ \frac{\partial^2}{\partial y^2}\Bigr) +iky\Bigl( \frac{\partial}{\partial x}+ i\frac{\partial}{\partial y}\Bigr)\]
splits as $\Delta_k = -\xi_{2-k}\xi_k$, we thus have $\Delta_{3/2}\mathcal{F}=0$ and so $\mathcal F$ is a harmonic Maass form.
Further, we have seen that the image under $\xi_{3/2}$, and hence the shadow of $q^{7/8}\,\mathcal{H}_{8,7}(q)$, is $-1/(4\pi \sqrt{2})\ko \oldtheta_2$.

Applying Lemma 3.7 in \cite{RZ}, one can immediately deduce that $\modtheta_{a,b}^{c_1,c_2}[f]| \gamma$ is holomorphic at $\infty$ for all $\gamma\in \operatorname{SL}_2(\Z)$. Since the order of the pole of $\eta^{-3}$ at $\infty$ is bounded, we get that $\mathcal{F}$ also satisfies the growth condition of a harmonic Maass form. 
\end{proof}



\begin{thebibliography}{99}


\bibitem{CG}
R.~Chen, F.~Garvan, \textit{A proof of the mod 4 unimodal sequence conjectures and related mock theta functions}, preprint (2020), arXiv:2010.14315.
  


\bibitem{GZ}
L.~G\"ottsche and D.~Zagier, \textit{Jacobi forms and the structure of Donaldson invariants for 4-manifolds with $b_+=1$}, Selecta Math.\ (N.S.) \textbf{4} (1998), no.~1, 69--115.



\bibitem{humbert}
G.~Humbert, \textit{Formules relatives aux nombres de classes des formes quadratiques binaires et positives}, Journ. de Math. \textbf{3} (1907), 337--449.



\bibitem{ogg}
A.~Ogg, \textit{Modular Forms and Dirichlet Series}, Mathematics lecture note series, W.A. Benjamin, 1969.

  
\bibitem{RZ}
C.~Roehrig and S.~Zwegers, \textit{Theta Series for Quadratic Forms of Signature $(n-1,1)$ with (Spherical) Polynomials}, preprint (2021), arXiv:2102.09329.
  
\bibitem{schoeneberg}
B.~Schoeneberg, \textit{Das Verhalten von mehrfachen Thetareihen bei Modulsubstitutionen}, Math.~Ann. \textbf{116} (1939), no.~1, 511--523.

\bibitem{shimura}
G.~Shimura, \textit{On modular forms of half integral weight}, Ann.~Math. \textbf{97} (1973), no.~3, 440--481.

\bibitem{siegel}
C.L.~Siegel, \textit{Indefinite quadratische {Formen} und {Funktionentheorie}. {I}.}, Math.~Ann. \textbf{124} (1951), 17--54.

\bibitem{vigneras1}
M.-F.~Vign\'eras, \textit{S\'eries th\^eta des formes quadratiques ind\'efinies}, S\'eminaire Delange-Pisot-Poitou. Th\'eorie des nombres \textbf{17} (1975--1976), no.~1, 1--3.

\bibitem{vigneras2}
\bysame, \textit{S\'eries th\^eta des formes quadratiques ind\'efinies}, In: Modular functions of one variable VI, Springer Lecture Notes \textbf{627} (1977), 227--239.

\bibitem{zagier2}
D.~Zagier, \textit{From quadratic functions to modular functions}, Number Theory in Progress: Proceedings of the International Conference on Number Theory organized by the Stefan Banach International Mathematical Center in Honor of the 60th Birthday of Andrzej Schinzel, Zakopane, Poland, June 30-July 9, 1997, K. Gy\"ory, H. Iwaniec and J. Urbanowicz (eds.) De Gruyter 2012, 1147--1178.

\bibitem{zwegers}
S.~Zwegers, \textit{Mock Theta Functions}, Ph.D. Dissertation (2002), Utrecht.
\end{thebibliography}
\end{document}